\documentclass[11pt,reqno]{amsart}
\usepackage{xifthen}
\usepackage{pkgfile}

\usepackage[comma,sort&compress, numbers]{natbib}

\allowdisplaybreaks

\title[Monochromatic Subgraphs in Dense Multiplex Networks]{Monochromatic Subgraphs in Randomly Colored Dense Multiplex Networks}

\author[Daros Andrade]{Mauricio Daros Andrade}
\address{Department of Statistics and Data Science, University of Pennsylvania, Philadelphia, USA} 
\email{daros@wharton.upenn.edu} 

\author[Bhattacharya]{Bhaswar B. Bhattacharya}
\address{Department of Statistics and Data Science, University of Pennsylvania, Philadelphia, USA} 
\email{bhaswar@wharton.upenn.edu} 

\begin{document}

\begin{abstract}
	Given a sequence of graphs $G_n$ and a fixed graph $H$, denote by $T(H, G_n)$ the number of monochromatic copies of the graph $H$ in a uniformly random $c$-coloring of the vertices of $G_n$. In this paper we study the joint distribution of a finite collection of monochromatic graph counts in networks with multiple layers (multiplex networks). Specifically, given a finite collection of graphs $H_1, H_2,  \ldots, H_d$ we derive the joint distribution of $(T(H_1, G_n^{(1)}), T(H_2, G_n^{(2)}), \ldots, T(H_d, G_n^{(d)}))$, where $\bm G_n = (G_n^{(1)}, G_n^{(2)}, \ldots, G_n^{(d)})$ is a collection of dense graphs on the same vertex set converging in the joint cut-metric. The limiting distribution is the sum of 2 independent components: a multivariate Gaussian and a sum of independent bivariate stochastic integrals. This extends previous results on the marginal convergence of monochromatic subgraphs in a sequence of graphs to the joint convergence of a finite collection of monochromatic subgraphs in a sequence of multiplex networks. Several applications and examples are discussed.
\end{abstract}



\maketitle

\section{Introduction}

Let $G_n = (V(G_n), E(G_n))$ be a sequence of graphs with vertex set $V(G_n)= [n]:= \{1,2,\cdots, n\}$ and edge set $E(G_n)$. Suppose the vertices of $G_n$ are colored uniformly at random with $c \geq 2$ colors, that is,
\begin{equation*}
	\P(v \in [n] \text{ is assigned color } a \in \{1, 2, \ldots, c\} ) = \frac{1}{c},
\end{equation*}
independently from the other vertices. Given such a coloring an edge $(u, v) \in E(G_n)$ is said to be {\it monochromatic} if both $u$ and $v$ are assigned the same color. Formally, if $X_v$ denotes the color of the vertex $v \in [n]$, then $(u, v) \in E(G_n)$ is monochromatic if and only if $X_u = X_v$.
%
Denote by $T(K_2,G_n)$ the number of monochromatic edges in $G_n$, that is,
\begin{equation}\label{eq:T_K2_G}
	T(K_2,G_n) :=\frac{1}{2}\sum_{1\leq u\neq v\leq n}a_{uv}(G_n)\I\{X_u=X_v\} ,
\end{equation}
where $A(G_n)= (a_{uv}(G_n))_{1 \leq u, v \leq n}$ is the adjacency matrix of $G_n$.
More generally, for a fixed connected graph $H = (V(H), E(H))$ with $V(H)=\{1,\cdots, |V(H)| \}$, denote by $T(H, G_n)$ the number of monochromatic copies of $H$ in $G_n$, that is,
\begin{equation}\label{def:T_H_G}
	T(H, G_n):= \frac{1}{|{\Aut}(H)|}\sum_{\bm{s}\in [n]_{|V(H)|}}\prod_{(a,b)\in E(H)}a_{s_as_b}(G_n)\I\{X_{=\bm{s}}\} ,
\end{equation}
where:
\begin{itemize}
	\item $[n]_{|V(H)|}$ is the set of all $|V(H)|$-tuples $\bm{s}=(s_1,\cdots, s_{|V(H)|})\in [n]^{|V(H)|}$ with distinct indices.\footnote{For a set $S$, the set $S^N$ denotes the $N$-fold cartesian product $S\times S\times \cdots \times S$.} Thus, the cardinality of $[n]_{|V(H)|}$ is $\frac{n!}{(n-|V(H)|)!}$.
	\item For any $\bm{s}=(s_1,\cdots, s_{|V(H)|})\in [n]_{|V(H)|}$,
	      \begin{equation*}
		      \I\{X_{=\bm{s}}\}:=\I\{X_{s_1}=\cdots = X_{s_{|V(H)|}}\}.
	      \end{equation*}
	\item ${\Aut}(H)$ is the \emph{automorphism group} of $H$, that is, the set of permutations $\sigma$ of the vertex set $V(H)$ such that $(x,y)\in E(H)$ if and only if $(\sigma(x),\sigma(y))\in E(H)$.
\end{itemize}
Note that
\begin{equation*}
	\E T(H, G_n) = \frac{N(H, G_n)}{c^{|V(H)| - 1}},
\end{equation*}
where $N(H, G_n)$ is the number of copies of $H$ in $G_n$.

Various problems in combinatorial probability, nonparametric statistics, and theoretical computer science involve the study of  \eqref{eq:T_K2_G} and \eqref{def:T_H_G}. The following are some examples:

\begin{example}[Birthday problem]
	Suppose $G_n$ is a friendship-graph (two people are connected by an edge in the graph if they are friends) which is colored uniformly randomly with $c=365$ colors (where the colors correspond to birthdays and the birthdays are assumed to be uniformly distributed across the year). In this case a monochromatic edge in $G_n$ corresponds to two friends with the same birthday. The birthday problem asks for the probability that there are two friends with the same birthday, that is, $\P(T(K_2, G_n)>0)$. This is a celebrated problem in elementary probability (see \cite{dasgupta2005matching,diaconis2002bayesian,mitzenmacher2017probability} and references therein), which appears in many different contexts, such as the study of coincidences \cite{diaconis1989methods}, graph coloring problems \cite{barbour1992poisson,birthdayexchangeability,BDM}, testing discrete distributions \cite{batu,paninski2008coincidence,diakonikolas2019collision}, security of hash functions \cite{bellare2004hash,stinson2005cryptography}, and the discrete logarithm problem \cite{bhattacharyadlp,galbraith2012non,discretelogarithm}, among others. A natural generalization of the birthday problem is to consider higher-order birthday matches (often referred to as multi-collisions), that is, the number of $r$-tuples of friends sharing the same birthday \cite{diaconis1989methods,nandi2007multicollision,generalizedbirthdayproblem,birthdaymulticollison}. This corresponds to studying the asymptotics of $T(K_r, G_n)$, the number of monochromatic $r$-cliques $K_r$ in the friendship graph $G_n$ \cite{bhattacharya2020second,BMmonochromatic}.
\end{example}

\begin{example}[Graph-based nonparametric 2-sample tests]
	One of the fundamental problems in statistical inference is to decide whether or not two data sets are generated from
	the same statistical model. This is the classical 2-sample problem which can be formally stated as follows: Given independent
	and identically distributed samples $\sX_{n} = \{X_1, X_2, \ldots, X_n\}$ and $\sY_{n} = \{Y_1, Y_2, \ldots, Y_n\}$ from two
	multivariate distributions $F$ and $G$, respectively, the two-sample problem is to test the hypothesis:
	\begin{equation}\label{eq:FG}
		H_0 : F = G \quad \text{versus} \quad H_1 : F \ne G.
	\end{equation}
	A theme that has emerged repeatedly in the development of nonparametric 2-sample tests is the use of random geometric graphs. This includes the celebrated Friedman-Rafsky test based on the minimum spanning tree (MST) \cite{fr}, tests based on nearest-neighbor (NN) graphs \cite{henzenn,schilling}, and optimal matchings \cite{rosenbaum}. Here, the idea is to construct a random geometric graph $\sG_n$ of  the pooled sample $\sX_{n}\cup \sY_{n}$, and reject the null hypothesis $H_0$ in \eqref{eq:FG} if the number of edges in $\sG_n$ with one end-point in sample 1 and another in sample 2 is `small'. In other words, $H_0$ is rejected when the number of {\it non-monochromatic} edges in $\sG_n$ is `small' (where the sample labels correspond to the colors), which is equivalent to rejecting $H_0$ when the number of monochromatic edges is `large'. Hence, understanding the null distribution of such a test statistic entails studying the asymptotic properties of $T(K_2, \sG_n)$ with $c=2$ colors.
\end{example}

\begin{example}[Quadratic Rademacher chaos] \label{example:quadratic}  When $c=2$, $T(K_2, G_n)$ centered by its mean can be expressed as a quadratic form in Rademacher variables as follows:
	\begin{align}\label{eq:K2quadratic}
		T(K_2,G_n) - \E T(K_2, G_n) & =\frac{1}{2}\sum_{1\leq u\neq v\leq n}a_{uv}(G_n) \left( \I\{X_u=X_v\} - \frac{1}{2} \right) \nonumber \\
		                            & = \frac{1}{4}\sum_{1\leq u\neq v\leq n}a_{uv}(G_n) (2\I\{X_u=1\} -1)(2\I\{X_v=1\} -1)  \nonumber \\
		                            & = \frac{1}{4}  \sum_{1\leq u\neq v\leq n} a_{uv}(G_n) Y_u Y_v ,
	\end{align}
	where $Y_u = 2\I\{X_u=1\} -1$, for $1 \leq u \leq n$, are i.i.d. Rademacher random variables. Hence, when $c=2$, $T(K_2, G_n) - \E T(K_2, G_n)$ is the {\it quadratic Rademacher chaos} \cite{limitrademacher,rademacher},
	which appears in various contexts, for example, in the Hamiltonian of the Ising model \cite{bmpotts,BM}, testing
	independence in autoregressive models \cite{beran_dependence}, and Ramsey theory \cite{anitconcentrationramsey}.
	A characterization of all distributional limits of \eqref{eq:K2quadratic}  has been obtained recently in \cite{bhattacharya2024fluctuations}.
\end{example}

A number of recent papers have studied the limiting distribution of $T(H, G_n)$ (appropriately centered and scaled), beginning with the work of \citet{BDM} which considered the case $H=K_2$. Specifically, \cite[Theorem 1.3]{BDM} shows that
$T(K_2, G_n)$ is asymptotically normal whenever the fourth moment of a suitably normalized version of $T(K_2, G_n)$ converges to 3 (the fourth moment of the standard normal distribution). This is an instance of the celebrated fourth-moment phenomenon, which was originally discovered in the seminal papers \cite{nualart2005central,nourdin2009stein} and has, since then, emerged as a unifying principle governing the asymptotic normality of random multilinear forms (see the book \cite{nourdin2012normal} for further details). Going beyond monochromatic edges, subsequently it was shown that the fourth-moment phenomenon continues to hold for  the asymptotic normality of $T(K_3, G_n)$ (the number of monochromatic triangles) when $c \geq 5$ \cite{bhattacharya2022normal}  and for $T(H, G_n)$, for any fixed graph $H$, when $c \geq 30$ \cite{subgraph2023fourth}. For $c=2$ a more refined result about the asymptotic normality of $T(H, G_n)$ has been obtained recently by \citet{subgraph2024characterizing} in terms of local influence-based conditions. While the aforementioned results provide sufficient conditions (some of which are also necessary) for the asymptotic normality of $T(H, G_n)$, there are many instances where $T(H, G_n)$ has a non-Gaussian limit. This, in particular, is often the case when $G_n$ is a sequence of dense graphs, that is, $|E(G_n)| = \Theta(n^2)$. In this case, to describe the limiting distribution of $T(H, G_n)$ it is convenient to adopt the framework of graph limit theory \cite{lovasz_book} and assume $G_n$ converges (in the cut-distance) to a graphon $W$. Then limiting distribution of $T(H, G_n)$ (suitably standardized) has both Gaussian and non-Gaussian components, where the non-Gaussian component is a (possibly) infinite weighted sum of independent centered chi-squared random variables with the weights determined by the spectral properties of a graphon derived from $W$ (see \cite[Theorem 1.4]{BDM} and \cite[Theorem 1.3]{BMmonochromatic}).

Given the above results a natural next step is to study the joint distribution of a finite collection of monochromatic subgraphs. Specifically, given a collection of $d \geq 1$ graphs $\cH  = (H_1, H_2,  \ldots, H_d)$, what is the limiting joint distribution of $(T(H_1, G_n), T(H_2, G_n), \ldots, T(H_d, G_n))$ (suitably standardized)?  More generally, one can consider the joint distribution of
\begin{equation*}
	\bm T(\cH, \bm G_n) := (T(H_1, G_n^{(1)}), T(H_2, G_n^{(2)}), \ldots, T(H_d, G_n^{(d)})),
\end{equation*}
where $\bm G_n = (G_n^{(1)}, G_n^{(2)}, \ldots, G_n^{(d)})$ is a collection of $d\geq 1$ graphs having the  same vertex set $[n]$. Multiple graphs sharing a common vertex set are known as multilayer networks (or multiplexes). The ubiquitous presence of complex relational data with many interdependencies has propelled the rapidly developing literature on multiplex networks (see \cite{bianconi2018multilayer,de2013mathematical,chandna2022edge,networksbiology,arroyo2021inference,lei2020consistent,lunagomez2021modeling,macdonald2022latent,kivela2014multilayer,skeja2024quantifying} and the references therein). In this paper we derive the asymptotic distribution of $T(\cH, \bm G_n)$ for a sequence of dense multiplexes $\bm{G}_n$. To describe the limiting distribution, we define a natural notion of convergence of multiplexes (in the joint cut-distance) and also invoke the framework of multiple stochastic integrals (for capturing the non-Gaussian dependencies between the different marginals). Under the assumption that the multiplex $\bm G_n$ converges in the joint cut-metric, we show that the limiting distribution of $T(\cH, \bm G_n)$ has 2 independent components: one of which is a multivariate Gaussian and the other is a sum of independent bivariate stochastic integrals (see Theorem \ref{thm:jointTHGn}). The stochastic integrals are with respect to the same underlying Brownian motion on $[0, 1]$, which captures the dependence between the different marginals. A key ingredient of the proof is an invariance principle which allows us to replace the coloring variables with appropriately chosen Gaussian variables, which might be of independent interest (see Section \ref{sec:invariance}). As a corollary of our general result, we obtain, by considering $G_n^{(1)} = G_n^{(2)} = \cdots = G_n^{(d)} = G_n$, the joint distribution of $(T(H_1, G_n), T(H_2, G_n), \ldots, T(H_d, G_n))$ whenever $G_n$ converges to a graphon. This extends the result in \cite{BMmonochromatic} from the marginal distribution of the number of monochromatic copies of a particular graph $H$, to the joint distribution of the number of monochromatic copies of a collection of graphs (see Corollary \ref{cor:jointTHGnW}). Another important special case is when $H_1= H_2, \ldots = H_d = K_2$ is the edge (see Corollary \ref{cor:K2joint}), which establishes the joint distribution of a finite collection of random quadratic forms in Rademacher variables (recall Example \ref{example:quadratic}). In Section \ref{sec:examples} we discuss various examples, which, in particular, includes the correlated Erd\H{o}s-R\'enyi multiplex that commonly appears in the study of graph matching problems (see Example \ref{example:correlatedrandomgraph}).

\subsection{Convergence of Graphs and Multiplexes}

The theory of graph limits \cite{graph_limits_I,graph_limits_II,lovasz_book} has received phenomenal attention in recent years. Here, we will recall the basic notions about the convergence of graph sequences.  For a detailed exposition of the theory of graph limits refer to Lov\' asz \cite{lovasz_book}.

For 2 graphs $F$ and $G$, define the homomorphism density of $F$ into $G$ by
\begin{equation*}
	t(F,G) :=\frac{|\hom(F,G)|}{|V (G)|^{|V (F)|}} ,
\end{equation*}
where $|\hom(F,G)|$ denotes the number of homomorphisms from $F$ to $G$. In fact, $t(F, G)$ is the proportion of maps $\phi: V (F) \rightarrow V (G)$ which define a graph homomorphism. To define the continuous analogue of graphs, consider $\sW$ to be the space of all measurable functions from $[0, 1]^2$ into $[0, 1]$ which are symmetric, that is, $W(x, y) = W(y,x)$, for all $x, y \in [0, 1]$. An element of $\sW$ is called a {\it graphon}. For a simple graph $F$ with $V (F)= \{1, 2, \ldots, |V(F)|\}$ and a graphon $W \in \sW$, define the homomorphism density of $F$ in $W$ by
\begin{equation*}
	t(F,W) =\int_{[0,1]^{|V(F)|}}\prod_{(i,j)\in E(F)}W(x_i,x_j) \mathrm dx_1\mathrm dx_2\cdots \mathrm dx_{|V(F)|} =  \E\lp\prod_{(i, j)\in E(F)}W(U_i, U_j)\rp ,
\end{equation*}
where $U_1, U_2, \ldots$ are i.i.d. $\mathrm{Unif}[0, 1]$.
It is easy to verify that $t(F, G) = t(F, W^{G})$, where $W^{G}$ is the {\it empirical graphon} associated with the graph $G$ which is defined as:
\begin{equation}\label{eq:emp_graph}
	W^G(x, y) :=  \I\{(\lceil |V(G)|x \rceil, \lceil |V(G)|y \rceil)\in E(G)\}.
\end{equation}
(In other words, to obtain the empirical graphon $W^G$ from the graph $G$, partition $[0, 1]^2$ into $|V(G)|^2$ squares of side length $1/|V(G)|$, and let $W^G(x, y)=1$ in the $(i, j)$-th square if $(i, j)\in E(G)$, and 0 otherwise.)

The basic definition of graph limit theory is in terms of the convergence of the homomorphism densities \cite{graph_limits_I,graph_limits_II,lovasz_book}. Specifically, a sequence of graphs $G_n$ {\it converges to a graphon $W$} if for every finite simple graph $F$,
\begin{equation}
	\lim_{n\rightarrow \infty}t(F, G_n) = t(F, W).
	\label{eq:graph_limit}
\end{equation}
This convergence can be captured through the cut-distance, which is defined as follows:

\begin{definition}\label{definition:Wconvergence}\cite{lovasz_book}
	The {\it cut-distance} between two graphons $W_1, W_2 \in \cW$ is defined as,
	\begin{equation}\label{eq:Wconvergence}
		||W_1-W_2||_{\square}:=\sup_{A, B \in [0, 1]}\left|\int_{A \times B} W_1(x, y)-W_2(x, y)\,  \mathrm dx \mathrm dy \right| ,
	\end{equation}
	for $W_1, W_2 \in \sW$.
\end{definition}

Using the cut distance, one can define an equivalence relation on $\sW$ as follows: $W_1 \sim W_2$ when $\inf_{\phi} \|W_1 - W_2^\phi\|_{\square} = 0$, where the infimum taken over all measure-preserving bijections $\phi: [0, 1] \rightarrow [0, 1]$ and $W_2^\phi(x, y):= W_2(\phi(x), \phi(y))$, for $x, y \in [0, 1]$.
The orbit of $W \in  \sW$ is the set of all functions $W^{\phi}$, as $\phi$ varies over all measure preserving bijections from $[0, 1] \rightarrow [0, 1]$. Denote by $\tilde{W}$ the closure of the orbit of $W$ in $(\sW, \|\cdot\|_\square)$. The quotient space $\tilde{\sW} = \{\tilde{W} : W \in \sW\}$ of closed equivalence classes is associated with the {\it cut-metric}:
\begin{equation*}
	\delta_\square(\tilde{W}_1, \tilde{W}_2 ):=\inf_{\phi} \| W_1 - W_2^{\phi}\|_{\square} .
\end{equation*}
The central result in graph limit theory is that a sequence of graphs $G_n$ converges to $W$ in the sense of \eqref{eq:graph_limit} if
and only if $\tilde{W}^{G_n}$ is a Cauchy sequence in the $\delta_{\square}$ metric. Moreover, a sequence of graphs $G_n$ converges to a graphon $W$ if and only if $\delta_{\square}(\tilde{W}^{G_n}, \tilde{W} ) \rightarrow 0$ (see \cite[Theorem 3.8]{graph_limits_I}).

In this paper we deal with sequences of multilayer networks  (multiplexes), where one has multiple graphs on the same set vertices. A multiplex with $d$ layers will be referred to as a $d$-multiplex. To define convergence of a sequence of $d$-multiplexes $\bm G_n$ with a common vertex set $V(\bm G_n)  = [n]$, we need the notion of joint cut-metric. Towards this, denote by $\sW^d$ the $d$-fold Cartesian product of $\sW$. An element of $\sW^d$ will be referred to as a {\it $d$-graphon}. We can define an equivalence relation on $\sW^d$ as follows: Given $\bm W=(W_1, W_2, \ldots, W_d) \in \sW^d$ and $\bm W'=(W_1', W_2', \ldots, W_d') \in \sW^d$, we say $\bm{W} \sim \bm W'$ if $\inf_{\phi} \sum_{j=1}^d \|  W_j^\phi - W_j'\|_{\square} = 0$, where, as before, the infimum is taken over all measure-preserving bijections $\phi: [0, 1] \rightarrow [0, 1]$. Denote by $\tilde{\bm{W}}$ the closure of the orbit of $\bm{W}$ and the quotient space of closed equivalence classes as $\tilde{\sW}^d = \{\tilde{\bm{W}} : \bm{W} \in \sW^d\}$. For $\tilde{\bm{W}}, \tilde{\bm{W}}' \in \tilde{\sW}^d$, define the {\it joint cut-metric} as
\begin{equation*}
	\delta_\square(\tilde{\bm{W}}, \tilde{\bm{W}}' ):=\inf_{\phi} \sum_{j=1}^d \| W_j^{\phi} - W_j'\|_{\square} .
\end{equation*}
We will say that a sequence of $d$-multiplexes $\bm G_n = (G_n^{(1)}, G_n^{(2)}, \ldots, G_n^{d})$ converges to the $d$-graphon $\bm{W} = (W_1, W_2, \ldots, W_d) \in \sW^d$ if
\begin{equation}\label{eq:distanceW}
	\delta_\square( \tilde{\bm{W}}^{\bm{G}_n}, \tilde{\bm{W}} ) \rightarrow 0,
\end{equation}
where $\bm{W}^{\bm{G}_n} = (W^{G_n^{(1)}}, \ldots, W^{G_n^{(d)}})$ is the empirical $d$-graphon associated with the $d$-multiplex $\bm{G}_n$.

\subsection{Statement of the Results} With the above preparations, we are now ready to state our main results. To this end, fix an integer $d \geq 1$ and a finite collection of connected graphs $\cH= (H_1, H_2, \ldots, H_d)$. Then for a sequence of  $d$-multiplexes $\bm G_n = (G_n^{(1)}, G_n^{(2)}, \ldots, G_n^{(d)})$, define the vector of standardized monochromatic copies as
\begin{align}\label{eq:GammaHGnd}
	\bm{\Gamma}(\cH, \bm{G}_n)
	= \begin{pmatrix}
		  \Gamma(H^{(1)},G_n^{(1)}) \\
		  \vdots                    \\
		  \Gamma(H^{(d)},G_n^{(d)})
	  \end{pmatrix} ,
\end{align}
where
\begin{equation}\label{eq:GammaHGn}
	\Gamma(H,G_n) := |{\Aut}(H)| c^{|V(H)|-\frac{3}{2}} \left\{  \frac{ T(H,G_n) - \E(T(H,G_n)) }{n^{|V(H)|-1} } \right\},
\end{equation}
for any graph $H=(V(H), E(H))$ and any graph sequence $G_n$ with $G_n=([n], E(G_n))$. Our aim is to derive the limiting distribution of $\bm{\Gamma}(\cH, \bm{G}_n)$ when $\bm G_n$ converges to a $d$-multiplex $\bm{W}$.
We begin by introducing the notion of the 2-point conditional kernel:

\begin{definition}[2-point conditional kernel]
	Given a finite graph $H = (V(H), E(H))$ and a graphon $W$, define the {\it 2-point conditional kernel of $H$ with respect to $W$} as follows:\footnote{Note that $W_H$ can take values greater than $1$; hence, technically, $W_H$ is not a graphon. However, $W_H \leq |V(H)|(|V(H)|-1)$, hence, $W_H \in \sW_1$, where $\sW_1$ is the space of bounded, symmetric, measurable functions from $[0, 1]^2$ to $[0, \infty)$. }
	\begin{equation}\label{eq:WH}
		W_H(x,y):=\sum_{1\leq a \neq b \leq |V(H)|}\E\lp \prod_{(i, j)\in E(H)} W(U_i, U_j) \;\middle|\; U_{a}=x, U_{b}=y\rp.
	\end{equation}
\end{definition}

To describe the asymptotic variance of $\bm{\Gamma}(\cH, \bm{G}_n)$ it is convenient to define the following graph operation.

\begin{figure}
	\centering
	\includegraphics[scale=0.75]{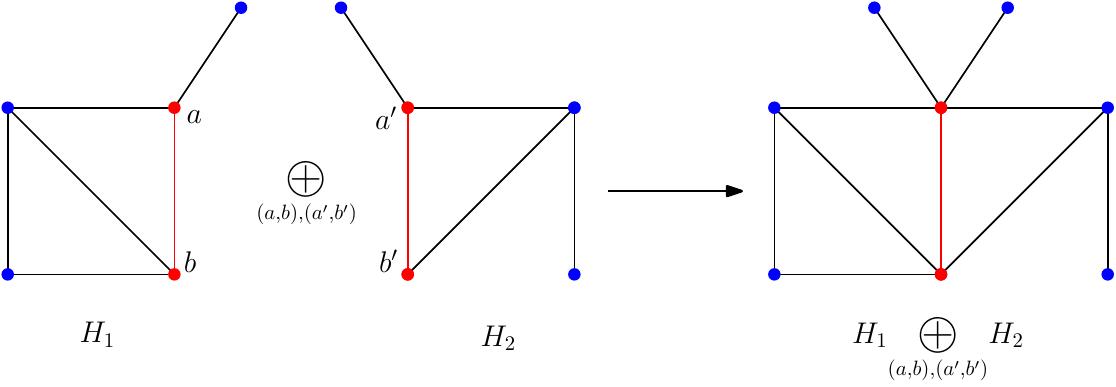}
	\caption{\small{ The $(a, b), (a', b')$-{\it join} of the graphs $H_1$ and $H_2$.} }
	\label{fig:vertexHab}
\end{figure}

\begin{definition}\label{definition:H2}
	Suppose $H_1 = (V(H_1), E(H_1))$ and $H_2 = (V(H_2), E(H_2))$ are two graphs with vertex sets $V(H_1) = \{1, 2, \ldots, |V(H_1)|\}$ and $V(H_2) = \{1, 2, \ldots, |V(H_2)|\}$ and edge sets $E(H_1)$ and $E(H_2)$, respectively. Then for $1 \leq a \ne b \leq |V(H_1)|$ and  $1 \leq a' \ne b' \leq |V(H_2)|$, the $(a, b), (a', b')$-{\it join} of $H_1$ and $H_2$ is the simple graph obtained by identifying the vertices $a$ and $b$ in $H_1$ with the vertices $a'$ and $b'$ in $H_2$, respectively. The resulting graph will be denoted by
	\begin{equation*}
		H_1 \bigotimes_{(a, b), (a', b')} H_2.
	\end{equation*}
	An example is shown in the Figure \ref{fig:vertexHab}. (Notice that the $(a,b),(a',b')$-join between 2 graphs can be different from the $(b,a),(a',b')$-join.)
\end{definition}

With the above notations we can now state our result about the asymptotic distribution of $\bm{\Gamma}(\cH, \bm{G}_n)$. Throughout, the stochastic integrals are considered in the Wiener-It\^{o} sense \cite{stochasticintegral} (see also \cite[Chapter 7]{SJIII}).

\begin{thm}\label{thm:jointTHGn}
	Suppose $\bm G_n$ is a sequence of $d$-multiplexes converging to a $d$-graphon $\bm{W}$. Let $\cH= (H_1, H_2, \ldots, H_d)$ be a collection of fixed graphs and $\bm{\Gamma}(\cH, \bm{G}_n)$ be as defined in \eqref{eq:GammaHGnd}. Assume, for $1  \leq i \ne j \leq d$, there exists $\rho_{ij} \geq 0$ such that
	\begin{equation}\label{eq:WGnHcovariance}
		\lim_{n \rightarrow \infty} \langle W_{H_i}^{G_n^{(i)}}, W_{H_j}^{G_n^{(j)}} \rangle = \rho_{ij} ,
	\end{equation}
	where $W_{H_i}^{G_n^{(i)}}$ is the 2-point conditional kernel of the empirical graphon $W^{G_n^{(i)}}$ with respect to $H_i$, for $1 \leq i \leq d$. Then, as $n \rightarrow \infty$,
	\begin{equation}\label{eq:distributionHGn}
		\bm{\Gamma}(\cH, \bm{G}_n)
		\xrightarrow{D}
		\sqrt{\frac{c-1}{2c}} \bm{Z} +\frac{1}{2\sqrt{c} }\sum_{a=1}^{c-1}\begin{pmatrix}\displaystyle\int_{[0,1]^2} (W_1)_{H_1} (x, y) \mathrm{d}B_x^{(a)}\mathrm{d}B_y^{(a)} \\ \vdots \\ \displaystyle\int_{[0,1]^2} (W_d)_{H_d}(x, y) \mathrm{d}B_x^{(a)}\mathrm{d}B_y^{(a)} \end{pmatrix} ,
	\end{equation}
	where
	\begin{itemize}
		\item $(W_i)_{H_i}$ is the 2-point conditional kernel of the graphon $W_i$ with respect to $H_i$,  for $1 \leq i \leq d$;
		\item $\bm{Z} \sim N_d(0, \Sigma)$, with $\Sigma = (\sigma_{ij})_{1 \leq i, j \leq d}$ defined as:
		      \everymath={\displaystyle}{
		      \begin{equation*}
			      \sigma_{ij} :=
			      \left\{
			      \begin{array}{cc}
				      \sum_{\substack{1\leq a \neq b\leq |V(H_i)| \\ 1\leq a' \ne b'\leq |V(H_i)|}}t\lp H_i\bigotimes_{(a, b), (a', b')}H_i, W_i \rp - \|(W_i)_{H_i}\|^2 &   \text{ for } i = j , \\
				      \rho_{ij} - \langle (W_i)_{H_i}, (W_j)_{H_j}\rangle & \text{ for } i \ne j ;
			      \end{array}
			      \right.
		      \end{equation*}
		      }
		\item $\{B^{(1)}_{t}, \ldots, B^{(c-1)}_{t} \}_{t \in [0, 1]}$ are i.i.d Brownian motions on $[0, 1]$ which are independent of $\bm{Z}$.
	\end{itemize}
\end{thm}

The proof of Theorem \ref{thm:jointTHGn} is given in Section \ref{sec:jointTHGnpf}. A key ingredient of the proof is a general invariance principle, which shows the following: random multilinear forms indexed by the variables which encode the vertex coloring are asymptotically close in moments to multilinear forms indexed by appropriately chosen Gaussian variables (see Section \ref{sec:invariance}). To apply this result, we first step express $T(H, G_n) - \E T(H, G_n)$, for any graph $H$, as a weighted sum of multilinear forms with degrees ranging from 2 to
$|V(H)|$ (see Lemma \ref{lm:T_H_G_n_multinomial_expansion}). Next, we show all terms in the expansion with degree greater than 2 are asymptotically negligible (see Lemma \ref{lem:Gamma_equal_Q_2}).
In other words, only the quadratic terms in the expansion of $T(H, G_n) - \E T(H, G_n)$ contribute to the limiting distribution. Then, using the invariance principle, the problem of deriving the joint distribution of $\bm{\Gamma}(\cH, \bm{G}_n)$ reduces to finding the joint distribution of certain Gaussian quadratic forms, which is analyzed using the Cram\'er-Wold device and spectral methods.

\begin{remark}
	Note that Theorem \ref{thm:jointTHGn} assumes that the layers of the multiplex $\bm{G}_n$ converge jointly in the cut-distance and the pairwise overlaps of the 2-point conditional kernels have a limit. In Section \ref{sec:examples} we provide examples which illustrate that both these assumptions are necessary for the limiting distribution of
	$\bm{\Gamma}(\cH, \bm{G}_n)$ to exist.
\end{remark}

One useful special case of Theorem \ref{thm:jointTHGn} is when $H_1= H_2 = \cdots = H_d = K_2$ is the edge. In this case, the 2-conditional kernel $W_{K_2}(x, y) = 2 W(x, y)$ (recall \eqref{eq:WH})
and, for $1 \leq i \ne j \leq d$,
\begin{equation*}
	\langle W_{K_2}^{G_n^{(i)}}, W_{K_2}^{G_n^{(j)}} \rangle = \frac{ 8 |E(G_n^{(i)}) \cap E(G_n^{(j)})|}{n^2} ,
\end{equation*}
where $E(G_n^{(i)}) $ is the edge set of the graph $G_n^{(i)}$, for $1 \leq i \leq d$. Hence, Theorem \ref{thm:jointTHGn} implies the following result:

\begin{cor}\label{cor:K2joint}
	Suppose $\bm G_n$ is a sequence of $d$-multiplexes converging to a $d$-graphon $\bm{W}$. Also, assume that there exists $\rho_{ij} \geq 0$ such that
	\begin{equation}\label{eq:covarianceedge}
		\frac{2 |E(G_n^{(i)}) \cap E(G_n^{(j)})|}{n^2} \rightarrow \rho_{ij} ,
	\end{equation}
	for $1 \leq i \ne j \leq d$. Then, as $n \rightarrow \infty$,
	\begin{equation*}
		\frac{1}{n}\lb
		\begin{pmatrix} T(K_2,G_{n}^{(1)})\\ \vdots \\ T(K_2,G_{n}^{(d)}) \end{pmatrix}
		- \frac{1}{c}\begin{pmatrix} |E(G_{n}^{(1)})|\\ \vdots\\ |E(G_{n}^{(d)})| \end{pmatrix}
		\rb
		\xrightarrow{D}
		\frac{\sqrt{c-1}}{\sqrt{2}c} \bm{Z} +\frac{1}{2c}\sum_{a=1}^{c-1}\begin{pmatrix}\displaystyle\int_{[0,1]^2}W_1\mathrm{d}B_x^{(a)}\mathrm{d}B_y^{(a)} \\ \vdots \\ \displaystyle\int_{[0,1]^2}W_d \mathrm{d}B_x^{(a)}\mathrm{d}B_y^{(a)} \end{pmatrix} ,
	\end{equation*}
	where
	\begin{itemize}
		\item $\bm{Z} \sim N_d(0, \Sigma)$, with $\Sigma = (\sigma_{ij})_{1 \leq i, j \leq d}$ defined as:
		      \everymath={\displaystyle}{
		      \begin{equation*}\sigma_{ij} :=
			      \left\{
			      \begin{array}{cc}
				      \int_{[0,1]^2} W_{i}(x, y)(1-W_{i}(x, y) ) \mathrm{d}x\mathrm{d}y     & \text{ for } i = j ,   \\
				      \rho_{ij} - \int_{[0,1]^2} W_i(x, y) W_j(x, y) \mathrm{d}x\mathrm{d}y & \text{ for } i \ne j ;
			      \end{array}
			      \right.
		      \end{equation*}
		\item $\{B^{(1)}_{t}, \ldots, B^{(c-1)}_{t} \}_{t \in [0, 1]}$ are i.i.d Brownian motions on $[0, 1]$ which are independent of $\bm{Z}$. }
	\end{itemize}
\end{cor}

If all the $d$ layers of $\bm{G}_n$ are the same, that is, $G_n^{(1)}= G_n^{(2)} = \cdots = G_n^{(d)} = G_n$, then Theorem \ref{thm:jointTHGn} gives the joint distribution of $(\Gamma(H_1, G_n), \ldots, \Gamma(H_d, G_n))$ for any collection of fixed graphs $\cH= (H_1, H_2, \ldots, H_d)$. In this case, the limit of $\langle W^{G_n}_{H_i}, W^{G_n}_{H_j} \rangle$, for $1 \leq i \ne j \leq d$, can be derived from the convergence of $G_n$ to a graphon $W$. Specifically, we have (see Lemma \ref{lm:WGnHvariance} for the proof),
\begin{equation*}
	\lim_{n \rightarrow \infty} \langle W^{G_n}_{H_1}, W^{G_n}_{H_2} \rangle = \sum_{\substack{1\leq u\neq v\leq |V(H_1)|\\ 1\leq u'\neq v'\leq |V(H_2)|}} t\lp H_1 \bigotimes_{(u,v),(u',v')} H_2, W \rp .
\end{equation*}
The above identity combined with Theorem \ref{thm:jointTHGn} gives the following result:

\begin{cor}\label{cor:jointTHGnW}
	Suppose $G_n$ is a sequence of graphs converging to a graphon $W$ and let $\cH= (H_1, H_2, \ldots, H_d)$ be a collection of fixed graphs. Then, as $n \rightarrow \infty$,
	\begin{equation*}
		\begin{pmatrix}
			\Gamma(H_1, G_n) \\
			\vdots           \\
			\Gamma(H_d, G_n)
		\end{pmatrix}
		\xrightarrow{D}
		\sqrt{\frac{c-1}{2c}} \bm{Z} +\frac{1}{2\sqrt{c} }\sum_{a=1}^{c-1}\begin{pmatrix}\displaystyle\int_{[0,1]^2} W_{H_1} (x, y) \mathrm{d}B_x^{(a)}\mathrm{d}B_y^{(a)} \\ \vdots \\ \displaystyle\int_{[0,1]^2} W_{H_d}(x, y) \mathrm{d}B_x^{(a)}\mathrm{d}B_y^{(a)} \end{pmatrix} ,
	\end{equation*}
	where
	\begin{itemize}
		\item $\bm{Z} \sim N_d(0, \Sigma)$, with $\Sigma = (\sigma_{ij})_{1 \leq i, j \leq d}$ defined as:
		      \everymath={\displaystyle}{
		      \begin{equation*}\sigma_{ij} :=
			      \left\{
			      \begin{array}{cc}
				      \sum_{\substack{1\leq a \neq b\leq |V(H_i)|    \\ 1\leq a' \ne b'\leq |V(H_i)|}}t\lp H_i\bigotimes_{(a, b), (a', b')}H_i, W \rp - \|W_{H_i}\|^2 &   \text{ for } i = j , \\
				      \  \sum_{\substack{1\leq a \neq b\leq |V(H_i)| \\ 1\leq a' \ne b'\leq |V(H_j)|}}t\lp H_i\bigotimes_{(a, b), (a', b')}H_j, W \rp - \langle W_{H_i}, W_{H_j}\rangle &   \text{ for } i \ne j ;
			      \end{array}
			      \right.
		      \end{equation*}
		      }
		\item $\{B^{(1)}_{t}, \ldots, B^{(c-1)}_{t} \}_{t \in [0, 1]}$ are i.i.d Brownian motions on $[0, 1]$ which are independent of $\bm{Z}$.
	\end{itemize}
\end{cor}

Note that when $\cH= (H)$ is a singleton we recover from Corollary \ref{cor:jointTHGnW} the marginal distribution of $T(H, G_n)$, which was proved in \cite[Theorem 3.1]{BMmonochromatic}.  In this case the limiting distribution can be alternately expressed in terms of the eigenvalues of the 2-point conditional kernel (recall \eqref{eq:WH}) as discussed in the following remark.

\begin{remark} Note that by the spectral theorem, for any graph $H$, the 2-point conditional kernel $W_H$ can be expressed as (see \cite[Section 7.5]{lovasz_book}):
	\begin{equation*}
		W_H (x, y) = \sum_{s=1}^\infty \lambda_s \phi_s(x) \phi_s(y),
	\end{equation*}
	where $\{\lambda_s\}_{s \geq 1}$ are the eigenvalues and $\{\phi_s\}_{s \geq 1}$ are an orthonormal collection of eigenvectors of the operator $\mathcal T_{W_H}$ defined by $\mathcal T_{W_H}(f)(x) = \int_0^1 W_H(x, y) f(y) \mathrm dy$. Hence, by the linearity of stochastic integrals and the product formula (see \cite[Page 100]{SJIII}), we have, for each $a \in [c]$,
	\begin{align*}
		\int_{[0, 1]^2} W_H(x, y) \mathrm d B_x^{(a)} \mathrm d B_y^{(a)} & = \sum_{s=1}^\infty \lambda_s \int_{[0, 1]^2} \phi_s(x) \phi_s(y) \mathrm d B_x^{(a)} \mathrm d B_y^{(a)}\nonumber \\
		                                                                  & = \sum_{s=1}^\infty \lambda_s \left( \left( \int_0^1 \phi_s(x) \mathrm d B_x^{(a)} \right)^2 - \int_0^1 \phi_s(x)^2 \mathrm d x \right) \nonumber \\
		                                                                  & = \sum_{s=1}^\infty \lambda_s \left( \left( \int_0^1 \phi_s(x) \mathrm d B_x^{(a)} \right)^2 - 1 \right) \tag*{ (by orthonormality) } \nonumber \\
		                                                                  & \stackrel{D} = \sum_{s=1}^\infty \lambda_s \left( Z_{s, a}^2 - 1 \right) ,
	\end{align*}
	where $Z_{s, a} \stackrel{D} =  \int_0^1 \phi_s(x) \mathrm d B_x$. Note that by orthonormality $\{Z_{s, a}\}_{s \geq 1, a \in [c]}$ is a collection of i.i.d. $N(0, 1)$ random variables. Hence, from Corollary \ref{cor:jointTHGnW} we obtain the following alternative expression of the limiting distribution of $\Gamma(H, G_n)$:
	\begin{align*}
		\Gamma(H, G_n) \dto
		\sqrt{\frac{c-1}{2c}} Z +\frac{1}{ 2 \sqrt{c} }\sum_{a=1}^{c-1} \sum_{s=1}^\infty \lambda_s \left( Z_{s, a}^2 - 1 \right) \nonumber \\
		\stackrel{D} =
		\sqrt{\frac{c-1}{2c}} Z +\frac{1}{ 2 \sqrt{c} } \sum_{s=1}^\infty \lambda_s \xi_s ,
	\end{align*}
	where $\{\xi_s\}_{ s \geq 1}$ is a collection of i.i.d. $\chi^2_{c-1} - (c-1)$ random variables that are independent of $Z \sim N(0, \sigma^2)$, with
	\begin{equation*}
		\sigma^2 =    \sum_{\substack{1\leq a \neq b\leq |V(H)|\\ 1\leq a' \ne b'\leq |V(H)|}}t\lp H\bigotimes_{(a, b), (a', b')}H, W \rp - \|W_{H}\|^2 .
	\end{equation*}
	This recovers the result in \cite[Theorem 1.3]{BMmonochromatic}.
\end{remark}

\subsection{Asymptotic Notation}
Throughout the paper we will use the following asymptotic notation: For two sequences $a_n$ and $b_n$ we say $a_n \lesssim b_n$ if for all $n$ large enough $a_n \leq C b_n$, for some constant $C > 0$. Similarly, $a_n \gtrsim b_n$ means  $a_n \geq C b_n$, for $n$ large enough and some constant $C > 0$. Subscripts in the above notation, for example $\lesssim_{\square}$ and $\gtrsim_{\square}$, denote that the hidden constants may depend on the subscripted parameters.


\section{A General Invariance Principle}
\label{sec:invariance}

For each $v \in [n]$ and $a\in [c]$, let
\begin{align}\label{eq:Xva}
	\tilde{X}_{v ,a} = \sqrt{c}\lp \I\{X_v=a\}-\frac{1}{c}\rp.
\end{align}
Let $\{Z_{v, a}\}_{v\in [n],a\in [c]}$ be a collection of i.i.d $N(0,1)$ and define
\begin{align}\label{eq:Zva}
	\tilde{Z}_{v,a}:= Z_{v,a} - \frac{1}{c}\sum_{a=1}^c Z_{v,a} .
\end{align}
We collect some easy facts about the random variables $\{\tilde{X}_{v, a}\}_{v\in [n],a\in [c]}$ and $\{\tilde{Z}_{v, a}\}_{v\in [n],a\in [c]}$ in the following observation:

\begin{obs}\label{obs:XZva}
	Let $\{\tilde{X}_{v, a}\}_{v\in [n],a\in [c]}$ and $\{\tilde{Z}_{v, a}\}_{v\in [n],a\in [c]}$ be as defined above. Then the following hold:

	\begin{enumerate}

		\item[$(a)$] \label{eq:first_second_moments_match:1} $\E \tilde{X}_{v,a} = \E \tilde{Z}_{v,a} =  0$, for $v\in [n],a\in [c]$.

		\item[$(b)$] \label{eq:first_second_moments_match:2} $\Var(\tilde{X}_{v,a}) = \Var(\tilde{Z}_{v,a}) = 1-\frac{1}{c}$, for $v\in [n],a\in [c]$.

		\item[$(c)$]  \label{eq:first_second_moments_match:3} $\Cov(\tilde{X}_{v,a},\tilde{X}_{v,b}) = \Cov(\tilde{Z}_{v,a},\tilde{Z}_{v,b}) = -\frac{1}{c}$, for $v\in [n]$ and $a \ne b \in [c]$.

		\item[$(d)$]  \label{eq:independenceuv} $\{\tilde{X}_{u,a}\}_{a\in [c]} \perp \{\tilde{X}_{v,a}\}_{a\in [c]}$ and $\{\tilde{Z}_{u,a}\}_{a\in [c]} \perp \{\tilde{Z}_{v,a}\}_{a\in [c]}$, for $u \ne v \in [n]$.

	\end{enumerate}
\end{obs}

Hereafter, denote
\begin{equation*}
	\bm{X}=(X_{v,a})_{v\in [n],a\in [c]}, ~ \tilde{\bm{X}}=(\tilde{X}_{v,a})_{v\in [n],a\in [c]}, ~ \bm{Z}=(Z_{v,a})_{v\in [n],a\in [c]} , \text{ and } \tilde{\bm{Z}}=(\tilde{Z}_{v,a})_{v\in [n],a\in [c]}.
\end{equation*}
Fix an integer $r\geq 0$. Then for a function $f:[n]_r \ra \R$, define
\begin{align}\label{eq:Tf}
	\cT(f; \tilde{\bm{X}}) := \frac{1}{n^{\frac{r}{2}}\sqrt{c}}\sum_{a=1}^{c}\sum_{\bm{s}\in [n]_{r}}f(\bm{s})\prod_{j=1}^{r} \tilde{X}_{s_j, a} .
\end{align}
Define $\cT(f; \bm{X})$ and $\cT(f; \tilde{\bm{Z}})$ similarly. Let $\cS_r$ denote the set of all permutations of the set $[r]$. Then define the symmetrization of $f$ as follows:
\begin{align*}
	\tilde{f} (s_1,\cdots s_r)=\frac{1}{r!}\sum_{\sigma \in \cS_r} f(s_{\sigma(1)},\cdots, s_{\sigma(r)}) .
\end{align*}
For functions $f,f':[n]_{r}\ra \R$, their inner product is defined as:
\begin{align*}
	\ip{f}{f'} := \frac{1}{n^r}\sum_{\bm{s}\in [n]_{r}}f(\bm{s})f'(\bm{s}) .
\end{align*}
and let $\|f\|:=\sqrt{\ip{f}{f}}$ be the associated norm. For a $\bm{s} \in [n]_{r}$, let $\underline{\bm{s}}$ denote the (unordered) set formed by the entries of $\bm s$. (For example, if $\bm s = (4, 2, 5)$, then $\underline{\bm{s}} = \{2, 4, 5\}$.) Notice that
\begin{align}\label{eq:ip_sym_S_s}
	\ip{\tilde{f}}{\tilde{f}'} = & \ip{\tilde{f}}{f'} = \ip{f}{\tilde{f}'} = \frac{1}{n^r}\frac{1}{r!}\sum_{\substack{\bm{s}, \bm{s}' \in [n]_{r} \\ \underline{\bm{s}}=\underline{\bm{s}}'}} f(\bm{s})f'(\bm{s}') .
\end{align}
We can then define the following bilinear operation:
\begin{align*}
	\ip{f}{f'}_{\bullet} := \ip{f}{\tilde{f}'} .
\end{align*}

\begin{lem}\label{lm:varianceTf}
	Fix integers $r, r' \geq 1$. Then for functions $f:[n]_{r}\ra \R$ and $f':[n]_{r'}\ra \R$, the following hold:
	\begin{itemize}
		\item If $r\neq r'$, then
		      \begin{align}\label{eq:cov_T_f_zero}
			      \Cov\lp \cT(f,\tilde{\bm{X}}),\cT(f',\tilde{\bm{X}}) \rp =0 .
		      \end{align}
		\item If $r=r'$, then
		      \begin{align*}
			      \Cov\lp \cT(f,\tilde{\bm{X}}),\cT(f',\tilde{\bm{X}}) \rp = \E \lp \cT(f,\tilde{\bm{X}}) \cT(f',\tilde{\bm{X}}) \rp
			      =\eta \cdot r!\ip{f}{f'}_{\bullet} ,
		      \end{align*}
		      where $\eta:= \lp1-\frac{1}{c}\rp^r+(c-1)\lp-\frac{1}{c}\rp^r$.
	\end{itemize}
	The same holds for $\Cov ( \cT(f,\tilde{\bm{Z}}),\cT(f',\tilde{\bm{Z}}) )$.
\end{lem}

\begin{proof}
	Note that for any $\bm{s} \in [n]_{r}$ and $a \in [c]$, $\E(\prod_{j = 1}^r \tilde{X}_{s_j, a})  = 0$. This implies, recall \eqref{eq:Tf}, $\E(\cT(f,\tilde{\bm{X}})) = 0$. Hence,
	\begin{align}\label{eq:covTf}
		\Cov\lp \cT(f,\tilde{\bm{X}}),\cT(f',\tilde{\bm{X}}) \rp & = \E \lp \cT(f,\tilde{\bm{X}}) \cT(f',\tilde{\bm{X}}) \rp \nonumber \\
		                                                         & =\frac{1}{n^{\frac{r+r'}{2}}c}\sum_{a,a'\in [c]} \sum_{\substack{\bm{s} \in [n]_{r} \bm{s}'\in [n]_{r'}}}f(\bm{s})f'(\bm{s}')\E\lp\prod_{j\in \underline{\bm{s}}}\tilde{X}_{s_j,a}\prod_{j\in \underline{\bm{s}}'}\tilde{X}_{s_j' ,a'} \rp ,
	\end{align}
	where $\bm{s} = (s_1, s_2, \ldots, s_r)$ and $\bm{s}' = (s_1', s_2', \ldots, s_{r'}')$, for $\bm{s} \in [n]_{r}$ and  $\bm{s}'\in [n]_{r'}$. Observe that if $\bm{s} \ne \bm{s'}$ (that is, either $\underline{\bm{s}} \backslash \underline{\bm{s}}'$ or $\underline{\bm{s}}' \backslash \underline{\bm{s}}$ is non-empty), then
	\begin{equation*}
		\E\lp\prod_{j\in \underline{\bm{s}}}\tilde{X}_{s_j,a}\prod_{j\in \underline{\bm{s}}'}\tilde{X}_{s_j' ,a'} \rp   = 0.
	\end{equation*}

	To begin with, suppose $r\neq r'$. Then $\bm{s} \ne \bm{s'}$, for any $\bm{s} \in [n]_{r}$ and  $\bm{s}'\in [n]_{r'}$. Hence, in this case all the terms in \eqref{eq:covTf} are zero, and the result in \eqref{eq:cov_T_f_zero} follows.

	Now, suppose $r=r'$ and $\bm{s}, \bm{s}' \in [n]_{r}$ is such that
	$\underline{\bm{s}}=\underline{\bm{s}}'$. Then
	\begin{align*}
		\frac{1}{c}\sum_{a,a'\in [c]}\E\lp\prod_{j\in \underline{\bm{s}}}\tilde{X}_{s_j, a}\prod_{j\in \underline{\bm{s}}'}\tilde{X}_{s_j, a'}\rp
		 & =\frac{1}{c}\sum_{a,a'\in [c]}\E\lp\prod_{j\in \underline{\bm{s}}}\tilde{X}_{s_j, a}\tilde{X}_{s_j, a'}\rp                                                                                               \\
		 & =\frac{1}{c}\sum_{a\in [c]}\E\lp\prod_{j\in \underline{\bm{s}}}\tilde{X}_{s_j, a}^2\rp + \frac{1}{c}\sum_{ a \ne a'\in [c] }\E\lp\prod_{j\in \underline{\bm{s}}}\tilde{X}_{s_j, a}\tilde{X}_{s_j, a'}\rp \\
		 & =\lp 1-\frac{1}{c}\rp^r +(c-1)\lp -\frac{1}{c}\rp^r = \eta .
	\end{align*}
	Hence, from \eqref{eq:covTf},
	\begin{align*}
		\Cov\lp \cT(f,\tilde{\bm{X}}),\cT(f',\tilde{\bm{X}}) \rp
		 & =\frac{1}{n^r c}\sum_{a,a'\in [c]} \sum_{\substack{ \bm{s} \in [n]_{r}, \bm{s}'\in [n]_{r} \\ \underline{\bm{s}}=\underline{\bm{s}}'}}f(\bm{s})f'(\bm{s}')\E\lp\prod_{j\in \underline{\bm{s}}}\tilde{X}_{s_j, a}\prod_{j\in \underline{\bm{s}}'}\tilde{X}_{s_j, a'} \rp\\
		 & =\eta \cdot \frac{1}{n^r}\sum_{\substack{ \bm{s} \in [n]_{r}, \bm{s}' \in [n]_{r'}         \\ \underline{\bm{s}}=\underline{\bm{s}}'}}f( \bm{s} ) f'( \bm{s}') \\
		 & =\eta \cdot r! \ip{f}{f'}_{\bullet} ,
	\end{align*}
	where the last step uses \eqref{eq:ip_sym_S_s}.
\end{proof}

Next, we show that $\cT(\cdot; \tilde{\bm{X}}) $ satisfies an invariance principle. Specifically, we show that the moments of $\cT(\cdot; \tilde{\bm{X}}) $ and $\cT(\cdot; \tilde{\bm{Z}}) $ are asymptotically close, for any finite collection of bounded functions. This result falls within the purview of universality/invariance principles for random multilinear forms (see, for example, \cite{chatterjee_invariance,stability,
	nourdin2010invariance,rotar2}).

\begin{ppn}\label{ppn:mthd_mmnts}
	Fix $K \geq 1$ and integers $r_1, r_2, \ldots, r_K \geq 0$. For $1 \leq k \leq K $, let $\{f_{n}^{(k)}\}_{n \geq 1}$ be sequence of functions such that $f_{n}^{(k)}:[n]_{r_k}\ra [-M, M]$, for some $M > 0$. Then
	\begin{align}\label{eq:mthd_mmnts_main_statement}
		\lim_{n\ra\infty}\left| \E\prod_{k=1}^K\cT\lp f_{n}^{(k)}; \tilde{\bm{X}}\rp - \E\prod_{k=1}^K\cT\lp f_{n}^{(k)}; \tilde{\bm{Z}}\rp \right| = 0 .
	\end{align}
\end{ppn}

\begin{proof}
	For notational convenience, denote
	$$\Delta_{n, K}:= \left| \E\prod_{k=1}^K\cT\lp f_{n}^{(k)}; \tilde{\bm{X}}\rp - \E\prod_{k=1}^K\cT\lp f_{n}^{(k)}; \tilde{\bm{Z}}\rp \right|.$$

	Note that
	\begin{align*}
		\Delta_{n, K} & =\frac{1}{n^{\frac{r_{\bullet}}{2}} c^{\frac{K}{2}}}\left| \sum_{\bm{a}\in [c]^K}\sum_{\substack{ \bm{s}^{(1)}, \ldots, \bm{s}^{(K)} \\ \bm{s}^{(k)}\in [n]_{r_k} }} \prod_{k = 1}^{K}f_{n}^{(k)}(\bm{s}^{(k)}) \left\{ \E \lp\prod_{k = 1}^{K}\prod_{j\in [n]_{r_k}}\tilde{X}_{s^{(k)}_j,a_k}\rp - \E\lp \prod_{k = 1}^{K}\prod_{j\in [n]_{r_k}}\tilde{Z}_{s^{(k)}_j,a_k}\rp \right\} \right| ,
	\end{align*}
	where  $r_{\bullet} = \sum_{k=1}^K r_k$, $\bm{a} = (a_1, a_2, \ldots, a_K) \in [c]^K$ and $\bm{s}^{(k)} = (s_1^{(k)}, s_2^{(k)}, \ldots, s_{r_k}^{(k)}) \in [n]_{r_K}$, for $1 \leq k \leq K$. Using the fact that $|f_n^{(k)}| \leq M$, for $1 \leq k \leq K$,  gives
	\begin{align*}
		\Delta_{n, K} & \leq \frac{M^K}{n^{\frac{r_{\bullet}}{2}}c^{\frac{K}{2}}} \sum_{\bm{a}\in [c]^K} \sum_{\substack{ \bm{s}^{(1)}, \ldots, \bm{s}^{(K)} \\ \bm{s}^{(k)}\in [n]_{r_k} }}  \left|\E \lp\prod_{k = 1}^{K}\prod_{j\in [n]_{r_k}}\tilde{X}_{s^{(k)}_j,a_k}\rp - \E\lp \prod_{k = 1}^{K}\prod_{j\in [n]_{r_k}}\tilde{Z}_{s^{(k)}_j,a_k}\rp \right|  .
		\stepcounter{equation}\tag{\theequation}
		\label{eq:mthd_mmnts_difference}
	\end{align*}
	Given $\bm{a} \in [c]^K$ and $\bm S : = (\bm{s}^{(1)},\cdots, \bm{s}^{(K)})$, such that $\bm{s}^{(k)} \in [n]_{r_k}$, for $1 \leq k \leq K$, define the edge-colored multi-hypergraph $F_{\bm{S}, \bm{a}} = (V(F_{\bm{S}, \bm{a}}), E(F_{\bm{S}, \bm{a}}))$ as follows:
	$$V(F_{\bm{S}, \bm{a}}) := \bigcup_{k=1}^K \underline{\bm{s}}^{(k)} \quad \text{ and } \quad E(F_{\bm{S_K}, \bm{a}} ) := \{\underline{\bm{s}}^{(1)}, \ldots, \underline{\bm{s}}^{(K)}\},$$
	and the color of the hyperedge $\underline{\bm{s}}^{(k)}$ is $a_k$, for $1 \leq k \leq K$. Note that the hyperedges in $\bm{s}^{(1)},\cdots, \bm{s}^{(K)}$ can be repeated, so $E(F_{\bm{S_K}, \bm{a}} )$ is a multi-set and $F_{\bm{S}, \bm{a}}$ is a multi-hypergraph. Also, define
	$$\tilde{\bm{X}}_{F_{\bm S, \bm a}} := \prod_{k = 1}^{K}\prod_{j\in [n]_{r_k}}\tilde{X}_{s^{(k)}_j,a_k} \text{ and } \tilde{\bm{Z}}_{F_{\bm S, \bm a}} := \prod_{k = 1}^{K}\prod_{j\in [n]_{r_k}}\tilde{Z}_{s^{(k)}_j,a_k} .$$
	The following lemma records 2 crucial properties of $\E \tilde{\bm{X}}_{F_{\bm S, \bm a}}$ and $\E \tilde{\bm{Z}}_{F_{\bm S, \bm a}}$.

	\begin{lem}\label{lm:XZ} Denote by $d_j$ the degree of vertex $j \in V(F_{\bm S, \bm a})$. Then the following hold:
		\begin{itemize}
			\item If there exists $j \in V(F_{\bm S, \bm a})$ such that
			      $d_j = 1$, then  $$\E \tilde{\bm{X}}_{F_{\bm S, \bm a}} =  \E \tilde{\bm{Z}}_{F_{\bm S, \bm a}} = 0.$$

			\item If $d_j = 2$, for all $j \in V(F_{\bm S, \bm a})$, then
			      $$\E \tilde{\bm{X}}_{F_{\bm S, \bm a}} =  \E \tilde{\bm{Z}}_{F_{\bm S, \bm a}} .$$
		\end{itemize}

	\end{lem}

	\begin{proof}  Note that
		\begin{align}\label{eq:XFs}
			\E \tilde{\bm{X}}_{F_{\bm S, \bm a}} = \E\lp \prod_{j\in V(F_{\bm S, \bm a})}\prod_{a=1}^c \tilde{X}_{j,a}^{\ell_{j,a}} \rp= \prod_{j\in V(F_{\bm S, \bm a})} \E\lp \prod_{a=1}^c \tilde{X}_{j,a}^{\ell_{j,a}} \rp,
		\end{align}
		where $\ell_{j,a}$ be the number of hyperedges with color $a \in [c]$ containing the vertex $j \in V(F_{\bm S, \bm a})$.
		Now, if $j \in V(F_{\bm S, \bm a})$ is such that $d_j = 1$, then there exists $b \in [c]$ such that $\ell_{j,b} =1$ and $\ell_{j,a} =0 $, for all $a \in [c]\backslash{b}$. This means, 
		$$\E\lp\prod_{a=1}^c \tilde{X}_{j,a}^{\ell_{j,a}} \rp  = \E \tilde{X}_{j,b} = 0$$ and, consequently, $\E \tilde{\bm{X}}_{F_{\bm S, \bm a}} = 0$.  By the same argument, we also have, $\E \tilde{\bm{Z}}_{F_{\bm S, \bm a}} = 0$.

		Next, suppose $d_j = 2$ for all $j \in V(F_{\bm S, \bm a})$. Then there are 2 possibilities:
		\begin{itemize}

			\item There exists $b \in [c]$ such that $\ell_{j,b} =2$ and $\ell_{j,a} =0 $, for all $a \in [c]\backslash\{b\}$. This means, 
			      $$\E\lp\prod_{a=1}^c \tilde{X}_{j,a}^{\ell_{j,a}} \rp  = \E \tilde{X}_{j,b}^2 = 1-\frac{1}{c} , $$
			      by Observation \ref{obs:XZva}.
			      By the same argument,
			      $$\E\lp\prod_{a=1}^c \tilde{Z}_{j,a}^{\ell_{j,a}} \rp  = \E \tilde{Z}_{j,b}^2 = 1-\frac{1}{c}.$$

			\item There exists $b \ne b' \in [c]$ such that $\ell_{j,b} = \ell_{j,b'} = 1$ and $\ell_{j,a} =0 $, for all $a \in [c]\backslash\{ b, b' \}$. This means,
			      $$\E\lp\prod_{a=1}^c \tilde{X}_{j,a}^{\ell_{j,a}} \rp  = \E \tilde{X}_{j,b} \tilde{X}_{j,b'} = -\frac{1}{c}, $$
			      by Observation \ref{obs:XZva}. By the same argument,
			      $$\E\lp\prod_{a=1}^c \tilde{Z}_{j,a}^{\ell_{j,a}} \rp  = \E \tilde{Z}_{j,b} \tilde{Z}_{j,b'} = -\frac{1}{c}.$$

		\end{itemize}
		Hence, if $d_j = 2$, then
		\begin{align*}
			\E\lp\prod_{a=1}^c \tilde{X}_{j,a}^{\ell_{j,a}} \rp = \E\lp\prod_{a=1}^c \tilde{Z}_{j,a}^{\ell_{j,a}} \rp .
		\end{align*}
		Hence, taking product over $j \in V(F_{\bm S, \bm a})$ (recall \eqref{eq:XFs}), gives $\E \tilde{\bm{X}}_{F_{\bm S, \bm a}} = \E \tilde{\bm{Z}}_{F_{\bm S, \bm a}} $.
	\end{proof}

	The next lemma gives bounds on $\E \tilde{\bm{X}}_{F_{\bm S, \bm a}}$ and $\E \tilde{\bm{Z}}_{F_{\bm S, \bm a}}$.

	\begin{lem}\label{lm:expectationXZbound} For any $\bm{a} = (a_1, a_2, \ldots, a_K) \in [c]^K$ and $\bm S = (\bm{s}^{(1)},\cdots, \bm{s}^{(K)})$, with $\bm{s}^{(k)} \in [n]_{r_k}$, for $1 \leq k \leq K$, we have
		\begin{align*}
			\E \tilde{\bm{X}}_{F_{\bm S, \bm a}} & \lesssim_{K, c} 1 .  \stepcounter{equation}\tag{\theequation}\label{eq:mthd_mmnts_product_of_vs_final}
		\end{align*}
		Moreover, $ \E \tilde{\bm{Z}}_{F_{\bm S, \bm a}} \lesssim_{r_1, \ldots, r_K} 1$.
	\end{lem}

	\begin{proof}
		To begin with, define $\cC_j : = \{a \in [c]: \ell_{j, a} \geq 1\}$, for $j \in V(F_{\bm S, \bm a})$. Then, for each $j \in V(F_{\bm S, \bm a})$,
		\begin{align}
			\left|\E\lp \prod_{a=1}^c \tilde{X}_{j,a}^{\ell_{j,a}} \rp\right|
			              & = \left|\E\lp \prod_{a \in \cC_j} \tilde{X}_{j,a}^{\ell_{j,a}} \rp\right| \nonumber                                                                                                      \\
			              & = c^{\frac{1}{2} \sum_{a \in \cC_j} \ell_{j, a}} \left| \E \prod_{a \in \cC_j} \lp \I \{X_{j}=a\}-\frac{1}{c}\rp^{\ell_{j, a}} \right| \nonumber                                         \\
			              & \leq c^{\frac{d_j}{2} }  \E \prod_{a \in \cC_j} \left|\I \{X_{j}=a\}-\frac{1}{c}\right|^{\ell_{j, a}} \nonumber                                                                          \\
			              & = c^{\frac{d_j}{2} } \sum_{b=1}^c \P(X_j = b) \lp 1-\frac{1}{c} \rp ^{\ell_{j, b}} \prod_{ a \in \cC_j \backslash \{b\}}\frac{1}{c^{\ell_{j, a}}} \nonumber                              \\
			              & \leq c^{\frac{d_j}{2}  - 1} \sum_{b=1}^c c^{- \sum_{a \in \cC_j \backslash \{b\} } \ell_{j, a}} \nonumber                                                                                \\
			\label{eq:TC} & = c^{\frac{d_j}{2}  - 1} \sum_{b \in \cC_j} c^{- \sum_{a \in \cC_j \backslash \{b\} } \ell_{j, a}} +  c^{\frac{d_j}{2}  - 1} \sum_{b \notin \cC_j} c^{- \sum_{a \in \cC_j } \ell_{j, a}} \\
			\label{eq:TK} & \leq  |\cC_j| c^{\frac{d_j}{2}  - |\cC_j|} +  c^{\frac{d_j}{2}  - |\cC_j|}  \lesssim_K  c^{\frac{d_j}{2}  - |\cC_j|} ,
		\end{align}
		where \eqref{eq:TC} uses $\ell_{j,a}\geq 1$, for $a\in \cC_j$, and \eqref{eq:TK} uses $|\cC_j| \leq K$, for $j \in V(F_{\bm{S}, \bm{a}})$. Hence, from \eqref{eq:XFs},
		\begin{align*}
			\left|\E \tilde{\bm{X}}_{\bm S, \bm a} \right| & \lesssim c^{\frac{1}{2}\sum_{j\in V(F_{\bm{S}, \bm{a}})} d_j -\sum_{j\in V(F_{\bm{S}, \bm{a}})}|\cC_j|} = c^{\frac{1}{2}\sum_{k\in [K]}r_k-\sum_{j\in V(F_{\bm{S}, \bm{a}})}(|\cC_j|-1)-|V(F_{\bm{S}, \bm{a}})|} \lesssim_{K, c} 1.
		\end{align*}
		This proves \eqref{eq:mthd_mmnts_product_of_vs_final}.

		Next, we bound $\E \tilde{\bm{Z}}_{F_{\bm S, \bm a}}$. First, observe that for each $j \in V(F_{\bm S, \bm a})$,
		\begin{align}
			\left|\E\lp \prod_{a=1}^c \tilde{Z}_{j,a}^{\ell_{j,a}} \rp\right|
			 & \leq  \E \prod_{a \in \cC_j} \left| Z_{j, a} - \frac{1}{c} \sum_{b=1}^c Z_{j, b} \right| ^{\ell_{j, a}}  \nonumber                                                                                \\
			 & \leq  \prod_{a \in \cC_j} \left( \E  \left| Z_{j, a} - \frac{1}{c} \sum_{b=1}^c Z_{j, b} \right|^{\ell_{j, a} |\cC_j | }  \right)^{\frac{1}{|\cC_j|}} \tag*{(by H\"older's inequality)} \nonumber \\
			 & \leq  \prod_{a \in \cC_j} \left( \left( 1-\frac{1}{c}  \right)^{\frac{\ell_{j, a} |\cC_j |}{2}} \E  \left|  N(0, 1) \right|^{\ell_{j, a} |\cC_j | }  \right)^{\frac{1}{|\cC_j|}} \nonumber        \\
			 & \leq  \prod_{a \in \cC_j} \left(  \E  \left|  N(0, 1) \right|^{\ell_{j, a} |\cC_j | }  \right)^{\frac{1}{|\cC_j|}}  \lesssim_K 1, \nonumber
		\end{align}
		since $\ell_{j, a} \leq K$ and $|\cC_j| \leq K$, for $j \in V(F_{\bm S, \bm a})$. This implies,
		\begin{align*}
			| \E \tilde{\bm{Z}}_{F_{\bm S, \bm a}} | = \left| \E\lp \prod_{j\in V(F_{\bm S, \bm a})}\prod_{a=1}^c \tilde{Z}_{j,a}^{\ell_{j,a}} \rp \right| = \prod_{j\in V(F_{\bm S, \bm a})} \left| \E\lp \prod_{a=1}^c \tilde{Z}_{j,a}^{\ell_{j,a}} \rp \right|  \lesssim_{r_1, r_2, \ldots, r_K} 1 ,
		\end{align*}
		since $|V(F_{\bm S, \bm a})| \leq \sum_{k=1}^K r_k = r_{\bullet}$. This completes the proof of Lemma \ref{lm:expectationXZbound}.
	\end{proof}

	Now, let $\cF$ be the set of colored hypergraphs with the following properties:
	\begin{itemize}
		\item $F \in \cF$ is a hypergraph with vertex set contained in $[r_{\bullet}]$ and $K$ hyperedges and every hyperedge of $F$ is assigned a color from $[K]$.

		\item all vertices of $F$ have degree at least two, and

		\item some vertex of $F$ has degree at least 3.

	\end{itemize}
	Note that for $F \in \cF$ there are $O(n^{|V(F)|})$ terms in the joint sum \eqref{eq:mthd_mmnts_difference} such that $F_{\bm S, \bm a}$ is isomorphic (as a colored hypergraph) to $F$. In such case, $\E\tilde{\bm{X}}_{F}=\E\tilde{\bm{X}}_{F_{\bm S, \bm a}}$ and $\E\tilde{\bm{Z}}_{F}=\E\tilde{\bm{Z}}_{F_{\bm S, \bm a}}$. This implies, from \eqref{eq:mthd_mmnts_difference} and Lemma \ref{lm:XZ},
	\begin{align}
		\Delta_{n, K} & \lesssim_{M, K, c} \frac{1}{n^{\frac{r_{\bullet}}{2}}} \sum_{F \in \cF} n^{|V(F)|} | \E \tilde{\bm{X}}_{F} - \E \tilde{\bm{Z}}_{F}  |\nonumber \\
		              & \lesssim_{M, K, c} \sum_{F \in \cF}n^{ |V(F)| -\frac{r_{\bullet}}{2}} .
		\label{eq:mthd_mmnts_sum_c_n_power}
	\end{align}
	Note that $r_{\bullet} = \sum_{k=1}^Kr_k=\sum_{e\in E(F)}|e|$. Also, for $F \in \cF$, we have that $d_j \geq 2$ for all $j\in V(F)$ with strict inequality for some $j$. Therefore,
	\begin{align*}
		\frac{1}{2} \sum_{k=1}^Kr_k = \frac{1}{2}\sum_{e\in E(F)}|e|=\frac{1}{2}\sum_{j\in V(F)} d_j > |V(F)|.
	\end{align*}
	Hence, each term in  \eqref{eq:mthd_mmnts_sum_c_n_power} is $o(1)$. Since $|\cF| \lesssim_{K, c} 1$, the result in Proposition \ref{ppn:mthd_mmnts} now follows.
\end{proof}

\section{Proof of Theorem \ref{thm:jointTHGn}}
\label{sec:jointTHGnpf}

Recall the definition of $T(H,G_n)$ from \eqref{def:T_H_G}, which is the number of monochromatic copies of $H$ in $G_n$:
\begin{align}\label{eq:recall_T_H_G_n}
	T(H,G_n)= \frac{1}{|{\Aut}(H)|}\sum_{\bm{s}\in [n]_{|V(H)|}}\prod_{(a,b)\in E(H)}a_{s_as_b}(G_n)\I\{X_{=\bm{s}}\} .
\end{align}
We begin by rewriting $T(H,G_n)$ in terms of $\cT(\cdot,\tilde{\bm{X}})$. To this end, we set $V(H) =  [|V(H)|]:= \{1, 2, \ldots, |V(H)|\}$. Also, for $J=\{i_1<\cdots < i_{|J|}\}\subseteq [|V(H)|]$ and $\bm{s} \in [n]_{|V(H)|}$, denote $\bm{s}_J = (s_{i_1},\cdots, s_{i_{|J|}})$.

\begin{lem}\label{lm:T_H_G_n_multinomial_expansion}
	For $J=\{i_1<\cdots < i_{|J|}\}\subseteq [|V(H)|]$, define $f_{H, J}^{G_n}:[n]_{|J|}\ra [0,1]$ as follows:
	\begin{align}\label{eq:fH}
		f_{H, J}^{G_n}(\bm{t}):=\frac{1}{n^{|V(H)|-|J|}}\sum_{\substack{\bm{s}\in [n]_{|V(H)|} \\ \bm{s}_J=\bm{t}}}\prod_{(a,b)\in E(H)}a_{\bm{s}_a \bm{s}_b}(G_n) .
	\end{align}
	Then
	\begin{align}\label{eq:expanded_T_H_G_n}
		T(H,G_n)-\E(T(H,G_n))= \frac{1}{|{\Aut}(H)|}\sum_{\substack{J\subseteq V(H) \\ |J|\geq 2}}\cT\lp f_{H, J}^{G_n}; \tilde{\bm{X}} \rp \frac{n^{|V(H)|-\frac{|J|}{2}}}{c^{|V(H)|-\frac{|J|}{2}-\frac{1}{2}}} .
	\end{align}
\end{lem}

\begin{proof}
	Notice that
	\begin{align*}
		\I\{X_{=\bm{s}}\} & = \sum_{a=1}^c \prod_{j=1}^{|V(H)|} \I\{X_{s_j}=a\}                                                             \\
		                  & = \sum_{a=1}^c \prod_{j=1}^{|V(H)|}\lp \I\{X_{s_j}=a\} -\frac{1}{c}+\frac{1}{c}\rp                              \\
		                  & = \sum_{a=1}^c \sum_{J\subseteq V(H)}\prod_{j\in J}\lp \I\{X_{s_j}=a\} -\frac{1}{c}\rp \frac{1}{c^{|V(H)|-|J|}} \\
		                  & = \sum_{a=1}^c \sum_{J\subseteq V(H)}\prod_{j\in J} \tilde{X}_{s_j, a}\frac{1}{c^{|V(H)|-\frac{|J|}{2}}}  .
		\stepcounter{equation}\tag{\theequation}\label{eq:indicator_x_equal_s}
	\end{align*}
	Then recalling \eqref{eq:Tf} we have,
	\begin{align}\label{eq:THfGn}
		\cT(f_{H, J}^{G_n}; \tilde{\bm{X}}) & = \frac{1}{n^{|V(H)| - \frac{|J|}{2}}\sqrt{c}}\sum_{a=1}^{c}\sum_{\bm{t}\in [n]_{|J|}} \sum_{\substack{\bm{s}\in [n]_{|V(H)|}                                                    \\ \bm{s}_J=\bm{t}}}\prod_{(a,b)\in E(H)}a_{\bm{s}_a \bm{s}_b}(G_n) \prod_{j \in J} \tilde{X}_{t_j, a} \nonumber \\
		                                    & = \frac{1}{n^{|V(H)| - \frac{|J|}{2}}\sqrt{c}}\sum_{a=1}^{c} \sum_{\bm{s}\in [n]_{|V(H)|} } \prod_{(a,b)\in E(H)}a_{\bm{s}_a \bm{s}_b}(G_n) \prod_{j \in J} \tilde{X}_{s_j, a} .
	\end{align}
	Hence,
	\begin{align}\label{eq:THfGnJsum}
		 & \frac{1}{|{\Aut}(H)|}\sum_{J\subseteq V(H)} \cT\lp f_{H, J}^{G_n}; \tilde{\bm{X}} \rp \frac{n^{|V(H)|-\frac{|J|}{2}}}{c^{|V(H)|-\frac{|J|}{2}-\frac{1}{2}}} \nonumber                                                       \\
		 & = \frac{1}{|{\Aut}(H)|}\sum_{J\subseteq V(H)} \sum_{a=1}^{c} \sum_{\bm{s}\in [n]_{|V(H)|} } \prod_{(a,b)\in E(H)}a_{\bm{s}_a \bm{s}_b}(G_n) \prod_{j \in J} \tilde{X}_{s_j, a} \frac{1}{c^{|V(H)|-\frac{|J|}{2}}} \nonumber \\
		 & = \frac{1}{|{\Aut}(H)|} \sum_{\bm{s}\in [n]_{|V(H)|} } \prod_{(a,b)\in E(H)}a_{\bm{s}_a \bm{s}_b}(G_n) \I\{X_{=\bm{s}}\} ,
	\end{align}
	where the last step uses \eqref{eq:indicator_x_equal_s}.

	Now, observe that, if $J=\{j_0\}$, for some $j_0 \in V(H)$,
	\begin{align*}
		\sum_{a=1}^c\prod_{j\in J}\lp \I\{X_{s_j}=a\} -\frac{1}{c}\rp=\sum_{a=1}^c\lp \I\{X_{s_{j_0}}=a\} -\frac{1}{c}\rp=0.
	\end{align*}
	This means (recall \eqref{eq:THfGn}) that
	\begin{equation}\label{eq:THfGnJone}
		\sum_{\substack{J\subseteq V(H)\\ |J|= 1}}\cT\lp f_{H, J}^{G_n}; \tilde{\bm{X}} \rp = 0.
	\end{equation}
	Also, notice that when $J = \emptyset$,
	\begin{equation}\label{eq:THfGnJempty}
		\frac{1}{|{\Aut}(H)|}  \cT\lp f_{H,\emptyset}^{G_n}; \tilde{\bm{X}} \rp \frac{n^{|V(H)|}}{c^{|V(H)|-\frac{1}{2}}} = \E(T(H,G_n)) .
	\end{equation}
	Combining \eqref{eq:THfGnJsum}, \eqref{eq:THfGnJone}, and \eqref{eq:THfGnJempty} and rearranging the sums the result in Lemma \ref{lm:T_H_G_n_multinomial_expansion} follows.
\end{proof}

Now, let $\{U_r\}_{r\in V(H)}$ be a collection of i.i.d uniform random variables on $[0,1]$. Define the counterpart of $W_H$ for finite $n$ as follows:

\begin{align*}
	\t{W}^{G_n}_{H}(x,y) & := \sum_{1\leq u\neq v\leq |V(H)|} \E\left(\I\{A_n\}\prod_{(i, j)\in E(H)}W^{G_n}(U_i, U_j)\, \middle|\, Z_u=x,Z_v=y\right) ,
	\stepcounter{equation}\tag{\theequation}
	\label{eq:WHGnxy}
\end{align*}
where $A_n =\{(\ceil{nU_r})_{r\in V(H)}\in [n]_{|V(H)|}\}$ is the event that the variables $\{U_r\}_{r\in V(H)}$ fall in distinct intervals when $[0, 1]$ is partitioned in a grid of size $1/n$. Then we define
\begin{align}
	Q_2(H, G_n) & := \frac{1}{ 2 n\sqrt{c}}\sum_{a=1}^c \sum_{1 \leq u \ne v \leq n} \t{W}^{G_n}_{H}\lp\frac{u}{n},\frac{v}{n}\rp  \tilde{Z}_{u, a} \tilde{Z}_{v, a} ,
	\stepcounter{equation}\tag{\theequation}\label{eq:QHGn}
\end{align}
where $\{\tilde{Z}_{v, a}\}_{v \in [n], a \in [c]}$ is defined in \eqref{eq:Zva}.

\begin{lem}\label{lem:Gamma_equal_Q_2}
	Let $\Gamma(H,G_n)$   and $Q_2(H, G_n)$ be as defined in \eqref{eq:GammaHGn} and \eqref{eq:QHGn}, respectively. Then
	\begin{align}\label{eq:Q2HGnW}
		Q_2(H, G_n) = \cT\lp \sum_{\substack{J\subseteq V(H) \\ |J|= 2}} f_{H, J}^{G_n}; \tilde{\bm{Z}} \rp .
	\end{align}
	Also, as $n \rightarrow \infty$,
	\begin{align}
		\Gamma(H,G_n) & = \cT\lp \sum_{\substack{J\subseteq V(H) \\ |J|= 2}} f_{H, J}^{G_n}; \tilde{\bm{X}} \rp + o_{L_2}(1) \stepcounter{equation}\tag{\theequation}\label{eq:Gamma_eq_Gamma_2} .
	\end{align}
\end{lem}

\begin{proof} Recalling \eqref{eq:fH} and \eqref{eq:WHGnxy} gives, for $1 \leq i \ne j \leq n$,
	\begin{align*}
		\sum_{\substack{J\subseteq V(H)              \\ |J|= 2}} f_{H, J}^{G_n}(i, j) =\frac{1}{n^{|V(H)|-2}}\sum_{1 \leq u < v\leq |V(H)|}\sum_{\substack{\bm{s}\in [n]_{|V(H)|}\\ s_u=i,\ s_v=j  } }\prod_{(a,b)\in E(H)}a_{\bm{s}_a \bm{s}_b}(G_n)
		 & = \frac{1}{2}\t{W}^{G_n}_{H}\lp x,y \rp ,
	\end{align*}
	for $x\in [\frac{i}{n},\frac{i+1}{n}]$ and $y\in [\frac{j}{n},\frac{j+1}{n}]$.  Therefore,
	\begin{align*}
		\cT\left( \sum_{\substack{J\subseteq V(H)                                                                                                                            \\ |J|= 2}} f_{H, J}^{G_n} ; \tilde{\bm{Z}} \right)
		 & = \frac{1}{ 2 n\sqrt{c}}\sum_{a=1}^c \sum_{1 \leq u \ne v \leq n} \t{W}^{G_n}_{H}\lp\frac{u}{n},\frac{v}{n}\rp  \tilde{Z}_{u, a} \tilde{Z}_{v, a} = Q_2(H, G_n) .
	\end{align*}
	This proves \eqref{eq:Q2HGnW}.

	Next, recalling \eqref{eq:GammaHGn} and from \eqref{eq:expanded_T_H_G_n} we have,
	\begin{align}\label{eq:GammaHGnJsum}
		\Gamma(H, G_n) & = \frac{c^{|V(H)|-\frac{3}{2}}}{n^{|V(H)|-1}} \left\{ T(H,G_n)-\E(T(H,G_n)) \right\} \nonumber \\
		               & = \frac{1}{|{\Aut}(H)|}\sum_{\substack{J\subseteq V(H)                                            \\ |J|\geq 2}}\cT\lp f_{H, J}^{G_n}; \tilde{\bm{X}} \rp \frac{c^{\frac{|J|}{2}-\frac{1}{2}} }{n^{\frac{|J|}{2}-1}} .
	\end{align}
	Note that, for any $J \subseteq V(H)$, with $|J| \geq 3$, $\E \cT(f_{H, J}^{G_n}; \tilde{\bm{X}}) = 0$ and, from Lemma \ref{lm:varianceTf},
	\begin{align*}
		\frac{1}{n^{|J| -2}} \mathrm{Var} (\cT(f_{H, J}^{G_n}; \tilde{\bm{X}}) ) & = \frac{\eta |J|!}{n^{|J| -2}} \|\tilde{f}_{H, J}^{G_n} \|_2^2  = o(1) ,
	\end{align*}
	since $0\leq \tilde{f}_{H, J}^{G_n} \leq 1$. Hence, $\cT(f_{H, J}^{G_n}) = o_{L_2}(1)$, for all $J \subseteq V(H)$ such that $|J| \geq 3$. This implies, from \eqref{eq:GammaHGnJsum},
	\begin{align*}
		\Gamma(H,G_n) & = \frac{c^{|V(H)|-\frac{3}{2}}}{n^{|V(H)|-1}} \sum_{\substack{J\subseteq V(H) \\ |J|= 2}}\cT\lp f_{H, J}^{G_n}; \tilde{\bm{X}} \rp \frac{n^{|V(H)|-\frac{|J|}{2}}}{c^{|V(H)|-\frac{|J|}{2}-\frac{1}{2}}} + o_{L_2}(1)\\
		              & = \cT\lp \sum_{\substack{J\subseteq V(H)                                      \\ |J|= 2}} f_{H, J}^{G_n}; \tilde{\bm{X}} \rp + o_{L_2}(1) . \stepcounter{equation}\tag{\theequation}\label{eq:Gamma_eq_Gamma_2}
	\end{align*}
	This proves \eqref{eq:Gamma_eq_Gamma_2} and completes the proof of Lemma \ref{lem:Gamma_equal_Q_2}.
\end{proof}

The following lemma represents $Q_2(H,G_n)$ as a sum of independent bivariate stochastic integrals.

\begin{lem}\label{lem:from_correlated_to_independent}
	Let $\{B^{(1)}(t), B^{(2)}(t), \ldots, B^{(c-1)}(t)\}_{t \in [0, 1]}$ be a collection of independent Brownian motions on $[0, 1]$ and $Q_2(H,G_n)$ be as defined in \eqref{eq:QHGn}. Then
	\begin{align*}
		Q_2(H,G_n) \stackrel{D} = \frac{1}{2\sqrt{c}}\sum_{a=1}^{c-1}\int_{[0, 1]^2} \t{W}^{G_n}_{H}(x,y)\mathrm{d}B^{(a)}_x\mathrm{d}B^{(a)}_y .
	\end{align*}
\end{lem}
\begin{proof} Let $\{Z_{v, a}\}_{v\in [n],a\in [c]}$ be a collection of i.i.d $N(0,1)$ random variables. For $1 \leq v \leq n$, let $\bm{Z}_v = (Z_{v, 1}, Z_{v, 2}, \ldots, Z_{v, c})^\top$.  Then, for $1 \leq v \leq n$, $\bm{M} \bm{Z}_v = \tilde{\bm{Z}}_v$, where $\tilde{\bm{Z}}_v = (\tilde{Z}_{v, 1}, \tilde{Z}_{v, 2}, \ldots, \tilde{Z}_{v, c})^\top$, $\bm M= \bm I-\frac{1}{c}\I^{\top}\I$, and $\I$ is the vector of ones of length $c$.
	Note that $\bm{M}^\top \bm{M} = \bm{M}$ and $\bm{M}$ has $c-1$ non-zero eigenvalues all of which are 1. Hence, we can write,
	\begin{align*}
		Q_2(H,G_n) & =\frac{1}{ 2 n\sqrt{c}}\sum_{a=1}^c\sum_{1\leq u \neq v \leq n}\t{W}^{G_n}_{H}\lp \frac{u}{n},\frac{v}{n}\rp\tilde{Z}_{u,a}\tilde{Z}_{v,a}    \\
		           & = \frac{1}{ 2 n\sqrt{c}}\sum_{1\leq u \neq v \leq n}\t{W}^{G_n}_{H}\lp\frac{u}{n},\frac{v}{n}\rp  \bm{Z}_u^\top \bm{M}^\top \bm{M} \bm{Z}_v   \\
		           & = \frac{1}{ 2 n\sqrt{c}}\sum_{1\leq u \neq v \leq n}\t{W}^{G_n}_{H}\lp\frac{u}{n},\frac{v}{n}\rp  \bm{Z}_u^\top \bm{M} \bm{Z}_v \nonumber     \\
		           & \stackrel{D}= \frac{1}{ 2 n\sqrt{c}}\sum_{1\leq u \neq v \leq n}\t{W}^{G_n}_{H}\lp \frac{u}{n},\frac{v}{n}\rp\sum_{a=1}^{c-1} Y_{u,a} Y_{v,a}
		\tag*{(where $\{Y_{v, a}\}_{v\in [n],a\in [c-1]}$ are i.i.d. $N(0, 1)$)}
		\nonumber                                                                                                                                                  \\
		           & \stackrel{D}{=} \frac{1}{ 2 \sqrt{c}}\sum_{a=1}^{c-1}\int_{[0, 1]^2} \t{W}^{G_n}_{H}(x,y)\mathrm{d}B^{(a)}_x\mathrm{d}B^{(a)}_y,
	\end{align*}
	where the last inequality follows from observing that $\frac{1}{n}Y_{u,a}Y_{v,a}\stackrel{D}{=}\int_{[0,1]^2}\I\{[\frac{u-1}{n},\frac{u}{n})\times [\frac{v-1}{n},\frac{v}{n})\} \mathrm dB_x \mathrm dB_y$ and that $\t{W}^{G_n}_{H}$ is a constant on $[\frac{u-1}{n},\frac{u}{n})\times [\frac{v-1}{n},\frac{v}{n})$, for all $1\leq u\neq v\leq n$.
\end{proof}

With the above preparations, we proceed with the proof of Theorem \ref{thm:jointTHGn}. For this, suppose $\cH= (H_1, H_2, \ldots, H_d)$ is a finite collection of graphs and $\bm G_n = (G_n^{(1)}, G_n^{(2)}, \ldots, G_n^{(d)})$ be a sequence of $d$-multiplexes converging in joint cut-norm to $\bm W = (W_1, W_2, \ldots, W_d)$. Recall the definition of $\bm{\Gamma}(\cH, \bm{G}_n)$ from \eqref{eq:GammaHGnd} and define $$\bm{Q}_2(\cH, \bm{G}_n) = (Q_2(H_1, G_n^{(1)}), Q_2(H_2,G_n^{(2)}), \ldots, Q_2(H_d, G_n^{(d)}))^\top.$$
Then, for any vector $\bm \alpha = (\alpha_1, \alpha_2, \ldots, \alpha_d) \in \R^d$, from Proposition \ref{ppn:mthd_mmnts} and Lemma \ref{lem:Gamma_equal_Q_2} we have, for all integers $r \geq 1$,
\begin{align*}
	\lim_{n \rightarrow \infty} \left| \E ((\bm \alpha^ \top \bm{\Gamma}(\cH, \bm{G}_n))^r) - \E (\bm \alpha^ \top  \bm{Q}_2(\cH, \bm{G}_n) )^r ) \right|  = 0.
\end{align*}
Therefore, by the Cram\'er-Wold device, to prove Theorem \ref{thm:jointTHGn} it suffices to derive the limiting distribution of $\bm \alpha^ \top \bm{Q}_2(\cH, \bm{G}_n)$.
For this, note that
\begin{align}\label{eq:sum_of_stoch_int_tilde}
	\bm{\alpha}^\top \bm{Q}_2(\cH, \bm{G}_n) = \sum_{i=1}^d \alpha_i Q_2(H_i, G_n^{(i)}) =  \frac{1}{ 2 \sqrt{c}}\sum_{a=1}^{c-1}I^{(a)}_2\lp \sum_{i=1}^d \alpha_i \t{W}^{G_n^{(i)}}_{H_i} \rp,
\end{align}
where $I^{(a)}_2(F):=\int_{[0,1]^2}F(x,y)\mathrm{d}B^{(a)}_x\mathrm{d}B^{(a)}_y$, for any bounded symmetric kernel $F \in \sW_1$.
Observe that for any $1 \leq i \leq d$,
$$\|\t{W}^{G_n^{(i)}}_{H_i}-W^{G_n^{(i)}}_{H_i}\|_2^2 \leq |V(H_i)|^2\P(A_n^c) . $$
Since $(\ceil{nU_i})_{i \in V(H)}$ is a uniform random vector over $[n]^{V(H)}$,
$$\P(A_n)=\frac{n(n-1)\cdots (n-|V(H)| +1 )}{n^{|V(H)|}}\ra 1.$$
Hence, \eqref{eq:sum_of_stoch_int_tilde} has the same limiting distribution as
\begin{align}\label{eq:sum_of_stoch_int}
	\frac{1}{ 2 \sqrt{c}} \sum_{a=1}^{c-1} I^{(a)}_2\lp \sum_{i=1}^d \alpha_i W^{G_n^{(i)}}_{H_i} \rp .
\end{align}
Now, let
\begin{align}\label{eq:WGnalpha}
	W_{\bm{\alpha}, n} := \sum_{i=1}^d \alpha_i W^{G_n^{(i)}}_{H_i} \text{ and }  W_{\bm \alpha} := \sum_{i=1}^d \alpha_i (W_i)_{H_i} .
\end{align}
Since $\bm{G}_n$ converges to $\bm{W}$ in joint cut-metric, there is a sequence of invertible measure preserving maps $\phi_n:[0,1]\ra [0,1]$ such that
$$\sum_{i=1}^d \|(W^{G_n^{(i)}})^{\phi_n}-W_i \|_{\square} \ra 0.$$
Define the map $\eta_{H_i} : W \rightarrow W_{H_i}$, for $1 \leq i \leq d$. Notice that $\eta_{H_i}((W^{G_n^{(i)}})^{\phi_n}) = \eta_{H_i}(W^{G_n^{(i)}})^{\phi_n}$.  Then by Lemma \ref{lem:W_H_lipschitz}, as $n \rightarrow \infty$,
\begin{align}\label{eq:Walphaconvergence}
	\| W_{\bm{\alpha}, n}^{\phi_n} - W_{\bm{\alpha}}\|_\square \rightarrow 0 .
\end{align}
Using \eqref{eq:WGnalpha} we can rewrite \eqref{eq:sum_of_stoch_int}  as
\begin{align}\label{eq:QHGnWlinear}
	\frac{1}{ 2 \sqrt{c}} \sum_{a=1}^{c-1} I^{(a)}_2\lp \sum_{i=1}^d \alpha_i W^{G_n^{(i)}}_{H_i} \rp & = \frac{1}{ 2 \sqrt{c}} \sum_{a=1}^{c-1} I^{(a)}_2\lp W_{\bm{\alpha}, n} \rp \nonumber                                                                         \\
	                                                                                                  & \stackrel{D} = \frac{1}{ 2 \sqrt{c}} \sum_{a=1}^{c-1} I^{(a)}_2\lp W_{\bm{\alpha}, n}^{\phi_n} \rp  \nonumber                                                  \\
	                                                                                                  & = \frac{1}{ 2 \sqrt{c}}\sum_{a=1}^{c-1}I_2^{(a)}(W_{\bm \alpha, n}^{\phi_n}-W_{\bm{\alpha}})+\frac{1}{ 2 \sqrt{c}}\sum_{a=1}^{c-1}I_2^{(a)}(W_{\bm{\alpha}}) .
\end{align}
By Lemma \ref{lem:almost_independent_L2} the two terms in the RHS above are asymptotically independent, hence, it suffices to derive the limiting distribution of the first term.

\begin{lem}\label{lm:I2clt} As $n \rightarrow \infty$,
	\begin{align}\label{eq:l2clt}
		\frac{1}{2 \sqrt{c}}\sum_{a=1}^{c-1}I_2^{(a)}(W_{\bm \alpha, n}^{\phi_n}-W_{\bm{\alpha}}) \dto \sqrt{\frac{c-1}{2c}} N(0, \bm \alpha^\top \Sigma \bm \alpha ) ,
	\end{align}
	where $\Sigma$ is defined in Theorem \ref{thm:jointTHGn}.
\end{lem}

\begin{proof}
	By the spectral theorem for bounded symmetric kernels (see \cite[Section 7.5]{lovasz_book}),
	$$\Delta_{\alpha, n}(x, y):= W_{\bm \alpha,n}^{\phi_n}(x,y)-W_{\bm{\alpha}}(x,y)=\sum_{s=1}^\infty \lambda^{(n)}_s \phi_s^{(n)}(x) \phi_s^{(n)}(y), $$
	where $\{\lambda^{(n)}_s\}_{s \geq 1}$ are the eigenvalues and $\{ \phi_s^{(n)}\}_{s \geq 1}$ are a set of orthonormal eigenvectors for the operator
	$$T_{\Delta_{\bm{\alpha}, n}} f (x) = \int_0^1 \Delta_{\bm{\alpha}, n}(x, y) f(y) \mathrm dy.$$
	Define, for $s \geq 1$,
	\begin{align*}
		\xi_s=\sum_{a=1}^{c-1} \left(\int_0^1 \phi_s^{(n)}(x)\mathrm{d}B^{(a)}_x \right)^2-(c-1).
	\end{align*}
	By orthonormality, the variables $\{\xi_s\}_{s \geq 1}$ have independent $\chi^2_{c-1}-(c-1)$ distributions. Hence,
	\begin{align*}
		\frac{1}{2\sqrt{c}}\sum_{a=1}^{c-1}I_2^{(a)}(W_{\bm \alpha,n}^{\phi_n}-W_{\bm{\alpha}}) \stackrel{D} = \frac{1}{2 \sqrt{c}}\sum_{s=1}^{\infty}\lambda^{(n)}_s \xi_s.
	\end{align*}
	The result in \eqref{eq:l2clt} then follows from Lemma \ref{lm:chisquareCLT}, if we show the following:
	\begin{align}\label{eq:lambdamaxsum}
		\lim_{n \rightarrow \infty} \max_{s \geq 1}|\lambda_s^{(n)}| = 0 \text{ and } \lim_{n \rightarrow \infty}  \sum_{s \geq 1}|\lambda_s^{(n)}|_2^2 = \bm{\alpha}^\top \Sigma \bm{\alpha}.
	\end{align}

	To show the first condition in \eqref{eq:lambdamaxsum}, we will use the identity $\sum_{s \geq 1}(\lambda^{(n)}_s )^4=t(C_4, \Delta_{n, \bm{\alpha}})$ (see \cite[Section 7.5]{lovasz_book}), where $C_4$ denotes the cycle of size $4$. This implies,
	$$\max_{s \geq 1}|\lambda_s^{(n)}|^4 \leq \sum_{s=1}^\infty (\lambda^{(n)}_s )^4=t(C_4, \Delta_{n, \bm{\alpha}}) \leq \|\Delta_{n, \bm{\alpha}} \|_\square \ra 0, $$
	by \eqref{eq:Walphaconvergence}.

	We now prove the second condition in \eqref{eq:lambdamaxsum}. For this, observe  that
	\begin{align*}
		\lim_{n\ra\infty}\|\Delta_{n, \bm{\alpha}}\|_2^2 & = \lim_{n\ra\infty}\|W_{\bm \alpha,n}^{\phi_n}-W_{\bm{\alpha}}\|_2^2                                                                          \\
		                                                 & =\lim_{n\ra\infty}\|W_{\bm \alpha, n} \|_2^2 - 2 \lim_{n\ra\infty} \langle W_{\bm \alpha,n}, W_{\bm{\alpha}} \rangle +\|W_{\bm{\alpha}}\|_2^2 \\
		                                                 & =\lim_{n\ra\infty}\|W_{\bm \alpha, n} \|_2^2-\|W_{\bm{\alpha}}\|_2^2,
	\end{align*}
	using the fact that $\lim_{n\ra\infty} \langle W_{\bm \alpha,n}, W_{\bm{\alpha}} \rangle = \|W_{\bm{\alpha}}\|_2^2$, which follows from \eqref{eq:Walphaconvergence} and by invoking \cite[Lemma 8.22]{lovasz_book}. Also,
	\begin{align*}
		\lim_{n\ra\infty}\|W_{\bm \alpha, n} \|_2^2 & = \lim_{n\ra\infty}\left\| \sum_{i=1}^d \alpha_i W^{G_n^{(i)}}_{H_i}  \right\|_2^2                                                           \\
		                                            & = \lim_{n\ra\infty} \left\{  \sum_{i=1}^d \alpha_i^2 \| W^{G_n^{(i)}}_{H_i} \|_2^2 +  \sum_{1 \leq i \ne j \leq d} \alpha_i \alpha_j \langle
		W^{G_n^{(i)}}_{H_i} , W^{G_n^{(j)}}_{H^{(j)}} \rangle  \right\}                                                                                                                            \\
		                                            & = \sum_{i=1}^d \alpha_i^2 \tau(H_i, W_i) +  \sum_{1 \leq i \ne j \leq d} \alpha_i \alpha_j \rho_{ij}   ,
	\end{align*}
	where $\tau(H_i, W_i)$ is as defined in Lemma \ref{lm:WGnHvariance} and $\rho_{ij}$ is as in \eqref{eq:WGnHcovariance}. Combining the above and observing that $$\|W_{\bm{\alpha}}\|_2^2 = \sum_{i=1}^d \alpha_i^2 \|(W_i)_{H_i}\|_2^2 + \sum_{1 \leq i \ne j \leq d} \alpha_i \alpha_j \langle (W_i)_{H_i}, W^{(j)}_{H^{(j)}} \rangle,$$ the second condition in \eqref{eq:lambdamaxsum} follows. This completes the proof of Lemma \ref{lm:I2clt}.
\end{proof}

The proof of Lemma \ref{lm:I2clt} uses the following result about the limiting value of $\|W^{G_n}_{H}\|_2^2$.

\begin{lem}\label{lm:WGnHvariance}
	For any fixed graph $H = (V(H), E(H))$ and a sequence of graphs $G_n$ converging to a graphon $W$,
	$$\lim_{n \rightarrow \infty} \|W^{G_n}_{H}\|_2^2 = \sum_{\substack{1\leq u\neq v\leq |V(H)|\\ 1\leq u'\neq v'\leq |V(H)|}} t\lp H\bigotimes_{(a, b),(a', b')} H, W \rp := \tau(H, W).$$
	Moreover, for any two fixed graphs $H_1 = (V(H_1), E(H_1))$ and $H_2 = (V(H_2), E(H_2))$ and a sequence of graphs $G_n$ converging to a graphon $W$,
	\begin{align}\label{lm:WGnHcovariance}
		\lim_{n \rightarrow \infty} \langle W^{G_n}_{H_1}, W^{G_n}_{H_2} \rangle = \sum_{\substack{1\leq u\neq v\leq |V(H_1)| \\ 1\leq u'\neq v'\leq |V(H_2)|}} t\lp H_1 \bigotimes_{(a, b),(a', b')} H_2, W \rp .
	\end{align}
\end{lem}

\begin{proof}
	Let $H'$ be an isomorphic copy of $H$. Then
	\begin{align*}
		\|W^{G_n}_{H}\|_2^2
		 & = \int_{[0,1]^2}\lb \sum_{1\leq a \neq b \leq |V(H)|} \E\lp \prod_{(i, j)\in E(H)} W^{G_n}(U_i,U_j) \;\middle|\; U_a=x,U_b=y  \rp \rb^2\mathrm{d}x \mathrm{d}y \\
		 & = \int_{[0,1]^2}\sum_{\substack{(a, b)\in V(H)_2                                                                                                               \\ (a',b')\in V(H')_2}}\E\lp \prod_{(i, j) \in  E(H)\cup E(H')} W^{G_n}(U_i,U_j) \;\middle|\; \begin{pmatrix} U_a \\ U_b\end{pmatrix} = \begin{pmatrix} U_{a'} \\ U_{b'}\end{pmatrix} = \begin{pmatrix} x \\ y\end{pmatrix}  \rp \mathrm{d}x \mathrm{d}y\\
		 & = \sum_{\substack{1\leq a \neq b \leq |V(H)|                                                                                                                   \\ 1\leq a' \neq b' \leq |V(H)|}} t\lp H\bigotimes_{(a,b),(a',b')} H, W^{G_n} \rp
		\rightarrow \tau(H, W) ,
	\end{align*}
	since $G_n$ converges to $W$. The result in \eqref{lm:WGnHcovariance} can be proved similarly.
\end{proof}

To complete the proof of Theorem \ref{thm:jointTHGn}, recall the representation in \eqref{eq:QHGnWlinear}. Then applying Lemma \ref{lm:I2clt} and Lemma \ref{lem:almost_independent_L2} shows that
\begin{align*}
	\frac{1}{2 \sqrt{c}} \sum_{a=1}^{c-1} I^{(a)}_2\lp \sum_{i=1}^d \alpha_i W^{G_n^{(i)}}_{H_i} \rp \stackrel{D} \rightarrow  \sqrt{\frac{c-1}{2c}} \bm{\alpha}^\top \bm{Z} + \frac{1}{2 \sqrt{c}}\sum_{a=1}^{c-1}I_2^{(a)}(W_{\bm{\alpha}}) ,
\end{align*}
where $\bm{Z}$ is as defined Theorem \ref{thm:jointTHGn} and the 2 terms in the RHS above are independent. Finally, recall (from \eqref{eq:sum_of_stoch_int_tilde} and \eqref{eq:sum_of_stoch_int}) that the LHS above has the same limiting distribution as $\bm{\alpha}^\top \bm{Q}_2(\cH, \bm{G}_n)$. Hence, by the Cram\'er-Wold device the result in Theorem \ref{thm:jointTHGn} follows. \hfill $\qed$

\section{Examples}
\label{sec:examples}

In this section we illustrate our results in various examples. A few of the examples will involve random multiplexes. To this end, it is worth noting that Theorems \ref{thm:jointTHGn} can be easily extended to random multiplexes, when the limits in \eqref{eq:distanceW} and \eqref{eq:WGnHcovariance} hold in probability, under the assumption that the multiplex and its coloring are jointly independent (see, for example, \cite[Lemma C.1]{BDM}).

\begin{example}\label{example:correlatedrandomgraph}
	A natural example of a random multiplex is the {\it correlated Erd\H{o}s-R\'enyi model}, which is a 2-multiplex where the edges are dependent across the different layers. This model emerged from the study of network privacy \cite{pedarsani2011privacy} and is the basic underlying model in graph matching problems (see \cite{mao2023random,mao2021testing} and
	the references therein). Specifically, correlated Erd\H{o}s-R\'enyi model $\bm{G}(n, p, q, \rho) = (G_n^{(1)}, G_n^{(2)})$ is a 2-multiplex with a common vertex set $[n]$, where independently for every $1 \leq i < j \leq n$, we have
	\begin{align*}
		\P( (i, j) \in E(G_n^{(1)})) = p, \quad \P( (i, j) \in E(G_n^{(2)})) = q
	\end{align*}
	and
	\begin{align*}
		\P( (i, j) \in E(G_n^{(1)}) \cap E(G_n^{(2)})) =  \rho + p q := p_{1, 2},
	\end{align*}
	for $p, q \in (0, 1)$ and $p_{1,2}  \in [\max\{p+q-1,0\}, \min\{p,q\}]$. In this case, $(G_n^{(1)}, G_n^{(2)})$ converges jointly to $(W_1, W_2)$, where $W_1 \equiv p$ and $W_2 \equiv q$. Therefore,
	for any 2 graphs $H_1$ and $H_2$,
	\begin{align}\label{eq:WHrandom}
		((W_1)_{H_1}, (W_2)_{H_2})=\lp |V(H_1)|  (|V(H_1)| -1) p^{|E(H_1)|},|V(H_2)| (|V(H_2)| -1 ) q^{|E(H_2)|}\rp.
	\end{align}
	To apply Theorem \ref{thm:jointTHGn}, it remains to find the limit (recall \eqref{eq:distributionHGn}):
	\begin{align*}
		\lim_{n \rightarrow \infty} \langle W_{H_1}^{G_n^{(1)}}, W_{H_2}^{G_n^{(2)}} \rangle,
	\end{align*}
	and to compute $\Sigma$.

	Note that, given $\bm{G}_n = (G_n^{(1)}, G_n^{(2)})$,
	\begin{align*}
		\langle W_{H_1}^{G_n^{(1)}}, W_{H_2}^{G_n^{(2)}} \rangle & = \int_{[0, 1]^2} W_{H_1}^{G_n^{(1)}} (x, y) W_{H_2}^{G_n^{(2)}} (x, y) \mathrm d x \mathrm d y \nonumber \\
		                                                         & = \sum_{\substack{(a, b)\in V(H_1)_2                                                                      \\ (a',b')\in V(H_2)_2}} B_{H_1, H_2}^{G_n^{(1)} G_n^{(2)}} ((a, b), (a', b'))  ,
	\end{align*}
	where
	\begin{align*}
		 & B_{H_1, H_2}^{G_n^{(1)} G_n^{(2)}} ((a, b), (a', b'))  \nonumber                                                                                                                                                     \\
		 & := \int_{[0, 1]^2} \bigg\{ \E\lp \prod_{(i, j)\in E(H_1)} W^{G_n^{(1)}}(U_i, U_j) \middle|\; \begin{pmatrix} U_a \\ U_b\end{pmatrix} = \begin{pmatrix} x \\ y\end{pmatrix} \rp \nonumber                             \\
		 & \hspace{1.5in} \E\lp \prod_{(i, j)\in E(H_2)} W^{G_n^{(2)}}(U_i,U_j) \;\middle|\; \begin{pmatrix} U_{a'} \\ U_{b'}\end{pmatrix} = \begin{pmatrix} x \\ y\end{pmatrix}  \rp \bigg\} \mathrm{d}x \mathrm{d}y \nonumber \\
		 & \stackrel{L^2} \rightarrow \begin{cases}
			                              p_{1,2} \cdot p^{|E(H_1)|-1}q^{|E(H_2)|-1} & \text{if } (a,b)\in E(H_1) \text{ and } (a',b')\in E(H_2) , \\
			                              p^{|E(H_1)|}q^{|E(H_2)|}                   & \text{ otherwise}.
		                              \end{cases}
	\end{align*}
	Therefore, from \eqref{eq:WHrandom},
	\begin{align}\label{eq:WH1H2random}
		\langle W_{H_1}^{G_n^{(1)}}, W_{H_2}^{G_n^{(2)}} \rangle  - \ip{(W_1)_{H_1}}{(W_2)_{H_2}}  \stackrel{L^2} \rightarrow |E(H_1)||E(H_2)| \rho \cdot p^{|E(H_1)|-1}q^{|E(H_2)|-1} .
	\end{align}
	This implies,
	\begin{align}\label{eq:sigmarandom}
		\Sigma = \begin{pmatrix}
			         |E(H_1)||E(H_2)| p^{|E(H_1)|+ E(H_2)|-1}(1-p)            & |E(H_1)||E(H_2)| \rho \cdot p^{|E(H_1)|-1}q^{|E(H_2)|-1} \\
			         |E(H_1)||E(H_2)| \rho \cdot p^{|E(H_1)|-1}q^{|E(H_2)|-1} & |E(H_1)||E(H_2)| q^{|E(H_1)|+E(H_2)|-1}  (1-q)\end{pmatrix} .
	\end{align}
	Combining \eqref{eq:WH1H2random} and \eqref{eq:sigmarandom} and applying Theorem \ref{thm:jointTHGn} gives,
	\begin{align*}
		\begingroup
		\renewcommand*{\arraystretch}{1.5}
		\begin{pmatrix}
			\Gamma(H_1,G_n^{(1)}) \\ \Gamma(H_2,G_n^{(2)})
		\end{pmatrix}
		\xrightarrow{D}
		\sqrt{\frac{c-1}{2c}} \begin{pmatrix} Z_1 \\ Z_2 \end{pmatrix}
		+\frac{1}{ 2 \sqrt{c}}
		\begin{pmatrix}
			|V(H_1)| (|V(H_1)| -1) p^{|E(H_1)|} \\  |V(H_2)| (|V(H_2)| -1) q^{|E(H_1)|}    \end{pmatrix}
		\endgroup
		\xi ,
	\end{align*}
	with $\xi \sim \chi^2_{c-1}-(c-1)$ which is independent of $(Z_1,Z_2) \sim N_2(\bm 0, \Sigma)$.
\end{example}

Next, we construct an example where the Gaussian component of the limiting distribution in \eqref{eq:distributionHGn} has a degenerate covariance matrix.

\begin{example}\label{example:correlated_erdos_renyi_triangle_complement}
	Suppose $G_n^{(1)} \sim G(n, p)$, where $p \in (0, 1)$,  and consider the multiplex $\bm G_n = (G_n^{(1)}, G_n^{(2)})$, where
	$G_n^{(2)}$ is the complement graph of $G_n^{(1)}$. Then, by pluging in $H_1=H_2=K_3$, $c=2$, $q=1-p$ and $p_{1,2}=0$ into Example \ref{example:correlatedrandomgraph} gives
	\begin{align*}
		\begin{pmatrix}
			\Gamma(K_3,G_n^{(1)}) \\ \Gamma(K_3, G_n^{(2)})
		\end{pmatrix}
		\xrightarrow{D}
		3\sqrt{p(1-p)}\begin{pmatrix} p^2Z \\ -(1-p)^2Z \end{pmatrix}
		+3\sqrt{2}\begin{pmatrix}
			          \displaystyle p^3 \\ \displaystyle (1-p)^3
		          \end{pmatrix} \xi ,
	\end{align*}
	with $\xi\sim \chi^2_1-1$ which is independent of $Z \sim N(0,1)$.
\end{example}

The next example shows that marginal convergence (in the cut-metric) of the graphs in the different layers of a multiplex $\bm{G}_n$ and the convergence of the pairwise overlaps \eqref{eq:WGnHcovariance} are not enough for the convergence of the joint distribution of $\Gamma(\cH, \bm{G}_n)$. One needs to assume that the graphs in the different layers converge jointly in the cut-metric (as in \eqref{eq:distanceW}) for the limiting distribution of $\Gamma(\cH, \bm{G}_n)$ to exist.

\begin{figure}
	\centering
	\includegraphics[scale=0.75]{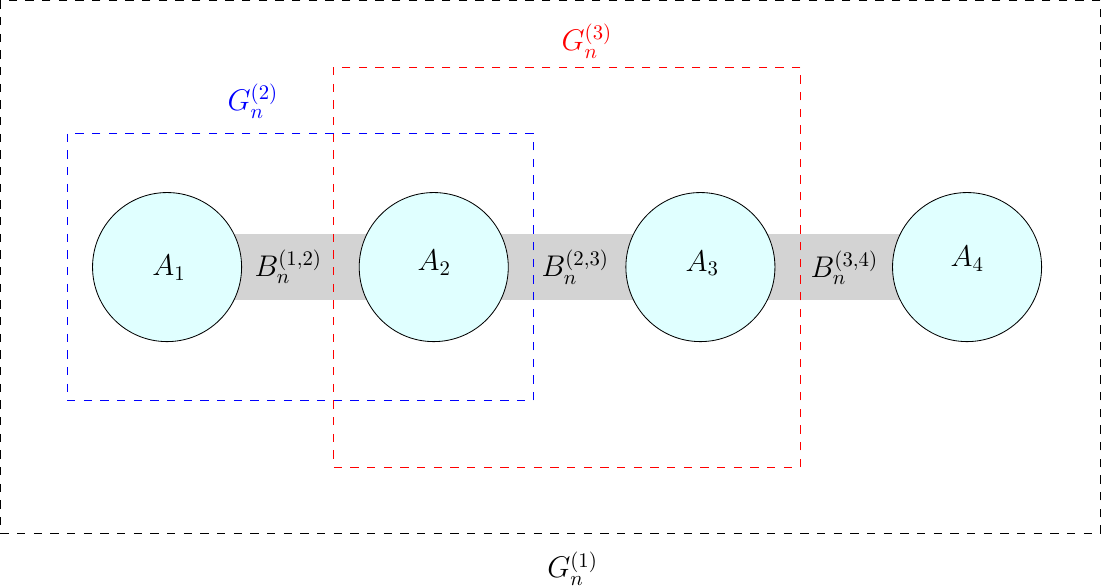}
	\caption{ \small{The graphs $G_n^{(1)}$,  $G_n^{(2)}$, and $G_n^{(3)}$ in Example \ref{ex:joint_convergence_necessary}. }}
	\label{fig:graphlayers}
\end{figure}

\begin{example}\label{ex:joint_convergence_necessary}
	Let $A_1, A_2, A_3, A_4$ be 4-disjoint sets of size $n$. Denote by $B_{n}^{(s, s+1)}$ the complete bipartite graph between the sets $A_s$ and $A_{s+1}$, for $1 \leq s \leq 3$ (see Figure \ref{fig:graphlayers}). Define\footnote{For 2 graphs $G_1 = (V(G_1), E(G_1))$ and $G_2 = (V(G_2), E(G_2))$, denote by $G_1 \bigcup G_2 = (V(G_1) \bigcup V(G_2), E(G_1) \bigcup E(G_2))$. }
	$$G_n^{(1)} = \bigcup_{s=1}^3 B_n^{(s, s+1)}, \quad G_n^{(2)} = B_n^{(1, 2)}, \text{ and } G_n^{(3)} = B_n^{(2, 3)}.$$
	(Note that, in other words, $G_n$ is the $n$-blow-up of a path with $4$ vertices (see \cite[Section 3.3]{lovasz_book}).) Now, consider the following 2 sequences of multiplexes:
	$$\bm G_n = (G_n^{(1)}, G_n^{(2)}) \text{ and } \tilde{\bm G}_n = (G_n^{(1)}, G_n^{(3)}).$$
	Note that the graph $G_n^{(1)}$ converges to the graphon $W_{1234}$ in Figure \ref{fig:bipartiteexample} (a) and the graphs $G_n^{(2)}$ and $G_n^{(3)}$ converge to the graphon $W_{23}$ in Figure \ref{fig:bipartiteexample} (b). Hence, marginally the graphs in the 2 layers of $\bm G_n$ and the graphs in the 2 layers of $\tilde{\bm G}_n$ converge to the same limits. Moreover,
	\begin{align*}
		\lim_{n \rightarrow \infty} \frac{2 |E(G_n^{(1)}) \cap E(G_n^{(2)})|}{n^2} = \lim_{n \rightarrow \infty}  \frac{2 |E(G_n^{(1)}) \cap E(G_n^{(3)})|}{n^2} = \frac{1}{8},
	\end{align*}
	that is, \eqref{eq:covarianceedge} converges to the same limit for both $\bm{G}_n$ and $\tilde{\bm{G}}_n$.

	\begin{figure}[H]
		\centering
		\begin{subfigure}[b]{0.45\textwidth}
			\centering
			\tikzset{every picture/.style={line width=0.75pt}} 

\begin{tikzpicture}[x=0.75pt,y=0.75pt,yscale=-1,xscale=1]

\draw  [color={rgb, 255:red, 153; green, 153; blue, 153 }  ,draw opacity=0.5 ] (15,15) -- (135,15) -- (135,135) -- (15,135) -- cycle ;
\draw    (15,150) -- (15,2) ;
\draw [shift={(15,0)}, rotate = 90] [color={rgb, 255:red, 0; green, 0; blue, 0 }  ][line width=0.75]    (10.93,-3.29) .. controls (6.95,-1.4) and (3.31,-0.3) .. (0,0) .. controls (3.31,0.3) and (6.95,1.4) .. (10.93,3.29)   ;
\draw    (0,135) -- (148,135) ;
\draw [shift={(150,135)}, rotate = 180] [color={rgb, 255:red, 0; green, 0; blue, 0 }  ][line width=0.75]    (10.93,-3.29) .. controls (6.95,-1.4) and (3.31,-0.3) .. (0,0) .. controls (3.31,0.3) and (6.95,1.4) .. (10.93,3.29)   ;
\draw  [color={rgb, 255:red, 0; green, 0; blue, 0 }  ,draw opacity=1 ][fill={rgb, 255:red, 0; green, 0; blue, 0 }  ,fill opacity=1 ][line width=0.75]  (45,105) -- (75,105) -- (75,135) -- (45,135) -- cycle ;
\draw  [color={rgb, 255:red, 0; green, 0; blue, 0 }  ,draw opacity=1 ][fill={rgb, 255:red, 0; green, 0; blue, 0 }  ,fill opacity=1 ][line width=0.75]  (15,75) -- (45,75) -- (45,105) -- (15,105) -- cycle ;
\draw  [color={rgb, 255:red, 0; green, 0; blue, 0 }  ,draw opacity=1 ][fill={rgb, 255:red, 0; green, 0; blue, 0 }  ,fill opacity=1 ][line width=0.75]  (45,45) -- (75,45) -- (75,75) -- (45,75) -- cycle ;
\draw  [color={rgb, 255:red, 0; green, 0; blue, 0 }  ,draw opacity=1 ][fill={rgb, 255:red, 0; green, 0; blue, 0 }  ,fill opacity=1 ][line width=0.75]  (75,75) -- (105,75) -- (105,105) -- (75,105) -- cycle ;
\draw  [color={rgb, 255:red, 0; green, 0; blue, 0 }  ,draw opacity=1 ][fill={rgb, 255:red, 0; green, 0; blue, 0 }  ,fill opacity=1 ][line width=0.75]  (105,45) -- (135,45) -- (135,75) -- (105,75) -- cycle ;
\draw  [color={rgb, 255:red, 0; green, 0; blue, 0 }  ,draw opacity=1 ][fill={rgb, 255:red, 0; green, 0; blue, 0 }  ,fill opacity=1 ][line width=0.75]  (75,15) -- (105,15) -- (105,45) -- (75,45) -- cycle ;

\end{tikzpicture} \\
			\small{(a)}
		\end{subfigure}
		\hfill
		\begin{subfigure}[b]{0.45\textwidth}
			\centering
			\tikzset{every picture/.style={line width=0.75pt}} 

\begin{tikzpicture}[x=0.75pt,y=0.75pt,yscale=-1,xscale=1]

\draw  [color={rgb, 255:red, 153; green, 153; blue, 113 }  ,draw opacity=0.5 ] (15,15) -- (135,15) -- (135,135) -- (15,135) -- cycle ;
\draw    (15,150) -- (15,2) ;
\draw [shift={(15,0)}, rotate = 90] [color={rgb, 255:red, 0; green, 0; blue, 0 }  ][line width=0.75]    (10.93,-3.29) .. controls (6.95,-1.4) and (3.31,-0.3) .. (0,0) .. controls (3.31,0.3) and (6.95,1.4) .. (10.93,3.29)   ;
\draw    (0,135) -- (148,135) ;
\draw [shift={(150,135)}, rotate = 180] [color={rgb, 255:red, 0; green, 0; blue, 0 }  ][line width=0.75]    (10.93,-3.29) .. controls (6.95,-1.4) and (3.31,-0.3) .. (0,0) .. controls (3.31,0.3) and (6.95,1.4) .. (10.93,3.29)   ;
\draw  [color={rgb, 255:red, 0; green, 0; blue, 0 }  ,draw opacity=1 ][fill={rgb, 255:red, 0; green, 0; blue, 0 }  ,fill opacity=1 ][line width=0.75]  (45,105) -- (75,105) -- (75,135) -- (45,135) -- cycle ;
\draw  [color={rgb, 255:red, 0; green, 0; blue, 0 }  ,draw opacity=1 ][fill={rgb, 255:red, 0; green, 0; blue, 0 }  ,fill opacity=1 ][line width=0.75]  (15,75) -- (45,75) -- (45,105) -- (15,105) -- cycle ;

\end{tikzpicture} \\
			\small{(b)}
		\end{subfigure}
		\hfill
		\caption{ \small{The graphons $W_{1234}$ and $W_{23}$ in Example \ref{ex:joint_convergence_necessary}.} }
		\label{fig:bipartiteexample}
	\end{figure}
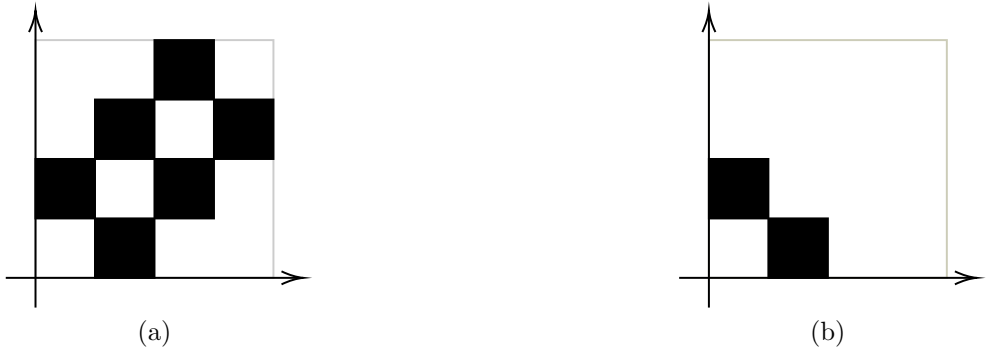
	Hence, the premises of Corollary \ref{cor:K2joint} hold for both $\bm G_n$ and $\tilde{\bm G}_n$. In particular, with $c=2$, we get
	\begin{align}\label{eq:Gn12}
		\begin{pmatrix} \Gamma(K_2,G_n^{(1)}) \\ \Gamma(K_2,G_n^{(2)}) \end{pmatrix}
		 & \xrightarrow{D}
		\begin{pmatrix}\displaystyle\int_{[0,1]^2}W_{1234}(x, y)\mathrm{d}B_x\mathrm{d}B_y \\ \displaystyle\int_{[0,1]^2}W_{23}(x, y)\mathrm{d}B_x\mathrm{d}B_y\end{pmatrix}
		\stackrel{D}{=}
		\frac{1}{4}\begin{pmatrix}\displaystyle Z_1Z_2+Z_2Z_3+Z_3Z_4 \\ Z_1Z_2 \end{pmatrix}
		=:
		\begin{pmatrix} Y \\ Y'\end{pmatrix}
	\end{align}
	and
	\begin{align}\label{eq:Gn13}
		\begin{pmatrix} \Gamma(K_2,G_n^{(1)}) \\ \Gamma(K_2,G_n^{(3)}) \end{pmatrix}
		 & \xrightarrow{D}
		\begin{pmatrix}\displaystyle\int_{[0,1]^2}W_{1234}(x, y)\mathrm{d}B_x\mathrm{d}B_y \\ \displaystyle\int_{[0,1]^2}W_{23}(x, y)\mathrm{d}B_x\mathrm{d}B_y\end{pmatrix}
		\stackrel{D}{=}
		\frac{1}{4}\begin{pmatrix}\displaystyle Z_1Z_2+Z_2Z_3+Z_3Z_4 \\ Z_2Z_3 \end{pmatrix}
		=:
		\begin{pmatrix} Y \\ Y'' \end{pmatrix} ,
	\end{align}
	where $Z_1, Z_2, Z_3, Z_4$ are i.i.d. $N(0,1)$ and $\{B_t\}_{t \in [0, 1]}$ is the standard Brownian Motion on $[0, 1]$. Note that the limiting distributions in \eqref{eq:Gn12} and \eqref{eq:Gn13} are different. In particular, $\E(Y-Y')^4=\frac{36}{4^4}$ which is different from $\E(Y-Y'')^4=\frac{24}{4^4}$. The reason for this is that in this case $\lim_{n\ra\infty}\delta_{\square}(\bm{W}^{\bm{G}_n},\bm{W}^{\t{\bm{G}}_n})\neq 0$.
\end{example}

\small

\subsection*{Acknowledgements}
B. B. Bhattacharya was supported by NSF CAREER grant DMS 2046393, NSF grant DMS 2113771, and a Sloan Research Fellowship. M. D. Andrade was supported by NSF grant 1953848.

\bibliographystyle{abbrvnat}
\bibliography{bibliography.bib}

\begin{thebibliography}{62}
\providecommand{\natexlab}[1]{#1}
\providecommand{\url}[1]{\texttt{#1}}
\expandafter\ifx\csname urlstyle\endcsname\relax
  \providecommand{\doi}[1]{doi: #1}\else
  \providecommand{\doi}{doi: \begingroup \urlstyle{rm}\Url}\fi

\bibitem[Arroyo et~al.(2021)Arroyo, Athreya, Cape, Chen, Priebe, and
  Vogelstein]{arroyo2021inference}
J.~Arroyo, A.~Athreya, J.~Cape, G.~Chen, C.~E. Priebe, and J.~T. Vogelstein.
\newblock Inference for multiple heterogeneous networks with a common invariant
  subspace.
\newblock \emph{Journal of Machine Learning Research}, 22\penalty0
  (142):\penalty0 1--49, 2021.

\bibitem[Barbour et~al.(1992)Barbour, Holst, and Janson]{barbour1992poisson}
A.~D. Barbour, L.~Holst, and S.~Janson.
\newblock \emph{{P}oisson approximation}, volume~2.
\newblock The Clarendon Press Oxford University Press, 1992.

\bibitem[Basak and Mukherjee(2017)]{bmpotts}
A.~Basak and S.~Mukherjee.
\newblock Universality of the mean-field for the {P}otts model.
\newblock \emph{Probability Theory and Related Fields}, 168\penalty0
  (3):\penalty0 557--600, 2017.

\bibitem[Batu et~al.(2013)Batu, Fortnow, Rubinfeld, Smith, and White]{batu}
T.~Batu, L.~Fortnow, R.~Rubinfeld, W.~D. Smith, and P.~White.
\newblock Testing closeness of discrete distributions.
\newblock \emph{Journal of the ACM (JACM)}, 60\penalty0 (1):\penalty0 1--25,
  2013.

\bibitem[Bellare and Kohno(2004)]{bellare2004hash}
M.~Bellare and T.~Kohno.
\newblock Hash function balance and its impact on birthday attacks.
\newblock In \emph{Advances in Cryptology-EUROCRYPT 2004: International
  Conference on the Theory and Applications of Cryptographic Techniques,
  Interlaken, Switzerland, May 2-6, 2004. Proceedings 23}, pages 401--418.
  Springer, 2004.

\bibitem[Beran(1972)]{beran_dependence}
R.~Beran.
\newblock Rank spectral processes and tests for serial dependence.
\newblock \emph{The Annals of Mathematical Statistics}, pages 1749--1766, 1972.

\bibitem[Bhattacharya(2016)]{bhattacharyadlp}
B.~B. Bhattacharya.
\newblock Collision times in multicolor urn models and sequential graph
  coloring with applications to discrete logarithms.
\newblock \emph{The Annals of Applied Probability}, 26\penalty0 (6):\penalty0
  3286--3318, 2016.

\bibitem[Bhattacharya and Mukherjee(2018)]{BM}
B.~B. Bhattacharya and S.~Mukherjee.
\newblock Inference in {I}sing models.
\newblock \emph{Bernoulli}, 24\penalty0 (1):\penalty0 493--525, 2018.

\bibitem[Bhattacharya and Mukherjee(2019)]{BMmonochromatic}
B.~B. Bhattacharya and S.~Mukherjee.
\newblock Monochromatic subgraphs in randomly colored graphons.
\newblock \emph{European Journal of Combinatorics}, 81:\penalty0 328--353,
  2019.

\bibitem[Bhattacharya et~al.(2017)Bhattacharya, Diaconis, and Mukherjee]{BDM}
B.~B. Bhattacharya, P.~Diaconis, and S.~Mukherjee.
\newblock Universal limit theorems in graph coloring problems with connections
  to extremal combinatorics.
\newblock \emph{The Annals of Applied Probability}, 27\penalty0 (1):\penalty0
  337, 2017.

\bibitem[Bhattacharya et~al.(2020)Bhattacharya, Mukherjee, and
  Mukherjee]{bhattacharya2020second}
B.~B. Bhattacharya, S.~Mukherjee, and S.~Mukherjee.
\newblock The second-moment phenomenon for monochromatic subgraphs.
\newblock \emph{SIAM Journal on Discrete Mathematics}, 34\penalty0
  (1):\penalty0 794--824, 2020.

\bibitem[Bhattacharya et~al.(2022)Bhattacharya, Fang, and
  Yan]{bhattacharya2022normal}
B.~B. Bhattacharya, X.~Fang, and H.~Yan.
\newblock Normal approximation and fourth moment theorems for monochromatic
  triangles.
\newblock \emph{Random Structures \& Algorithms}, 60\penalty0 (1):\penalty0
  25--53, 2022.

\bibitem[Bhattacharya et~al.(2024)Bhattacharya, Das, Mukherjee, and
  Mukherjee]{bhattacharya2024fluctuations}
B.~B. Bhattacharya, S.~Das, S.~Mukherjee, and S.~Mukherjee.
\newblock Fluctuations of quadratic chaos.
\newblock \emph{Communications in Mathematical Physics}, 405\penalty0
  (10):\penalty0 237, 2024.

\bibitem[Bianconi(2018)]{bianconi2018multilayer}
G.~Bianconi.
\newblock \emph{Multilayer networks: structure and function}.
\newblock Oxford university press, 2018.

\bibitem[Blei and Janson(2004)]{limitrademacher}
R.~Blei and S.~Janson.
\newblock Rademacher chaos: tail estimates versus limit theorems.
\newblock \emph{Arkiv f{\"o}r Matematik}, 42\penalty0 (1):\penalty0 13--29,
  2004.

\bibitem[Borgs et~al.(2008)Borgs, Chayes, Lov{\'a}sz, S{\'o}s, and
  Vesztergombi]{graph_limits_I}
C.~Borgs, J.~T. Chayes, L.~Lov{\'a}sz, V.~T. S{\'o}s, and K.~Vesztergombi.
\newblock Convergent sequences of dense graphs {I}: Subgraph frequencies,
  metric properties and testing.
\newblock \emph{Advances in Mathematics}, 219\penalty0 (6):\penalty0
  1801--1851, 2008.

\bibitem[Borgs et~al.(2012)Borgs, Chayes, Lov{\'a}sz, S{\'o}s, and
  Vesztergombi]{graph_limits_II}
C.~Borgs, J.~T. Chayes, L.~Lov{\'a}sz, V.~T. S{\'o}s, and K.~Vesztergombi.
\newblock Convergent sequences of dense graphs {II}. {M}ultiway cuts and
  statistical physics.
\newblock \emph{Annals of Mathematics}, 176\penalty0 (1):\penalty0 151--219,
  2012.

\bibitem[Cerquetti and Fortini(2006)]{birthdayexchangeability}
A.~Cerquetti and S.~Fortini.
\newblock A {P}oisson approximation for coloured graphs under exchangeability.
\newblock \emph{Sankhy{\=a}: The Indian Journal of Statistics}, pages 183--197,
  2006.

\bibitem[Chandna et~al.(2022)Chandna, Janson, and Olhede]{chandna2022edge}
S.~Chandna, S.~Janson, and S.~C. Olhede.
\newblock Edge coherence in multiplex networks.
\newblock \emph{arXiv:2202.09326}, 2022.

\bibitem[Chatterjee(2005)]{chatterjee_invariance}
S.~Chatterjee.
\newblock A simple invariance theorem.
\newblock \emph{arXiv preprint math/0508213}, 2005.

\bibitem[Das et~al.(2023)Das, Himwich, and Mani]{subgraph2023fourth}
S.~Das, Z.~Himwich, and N.~Mani.
\newblock A fourth-moment phenomenon for asymptotic normality of monochromatic
  subgraphs.
\newblock \emph{Random Structures \& Algorithms}, 63\penalty0 (4):\penalty0
  968--996, 2023.

\bibitem[DasGupta(2005)]{dasgupta2005matching}
A.~DasGupta.
\newblock The matching, birthday and the strong birthday problem: a
  contemporary review.
\newblock \emph{Journal of Statistical Planning and Inference}, 130\penalty0
  (1-2):\penalty0 377--389, 2005.

\bibitem[De~Domenico et~al.(2013)De~Domenico, Sol{\'e}-Ribalta, Cozzo,
  Kivel{\"a}, Moreno, Porter, G{\'o}mez, and Arenas]{de2013mathematical}
M.~De~Domenico, A.~Sol{\'e}-Ribalta, E.~Cozzo, M.~Kivel{\"a}, Y.~Moreno, M.~A.
  Porter, S.~G{\'o}mez, and A.~Arenas.
\newblock Mathematical formulation of multilayer networks.
\newblock \emph{Physical Review X}, 3\penalty0 (4):\penalty0 041022, 2013.

\bibitem[Diaconis and Holmes(2002)]{diaconis2002bayesian}
P.~Diaconis and S.~Holmes.
\newblock A {B}ayesian peek into {F}eller volume {I}.
\newblock \emph{Sankhy{\=a}: The Indian Journal of Statistics, Series A}, pages
  820--841, 2002.

\bibitem[Diaconis and Mosteller(1989)]{diaconis1989methods}
P.~Diaconis and F.~Mosteller.
\newblock Methods for studying coincidences.
\newblock \emph{Journal of the American Statistical Association}, 84\penalty0
  (408):\penalty0 853--861, 1989.

\bibitem[Diakonikolas et~al.(2019)Diakonikolas, Gouleakis, Peebles, and
  Price]{diakonikolas2019collision}
I.~Diakonikolas, T.~Gouleakis, J.~Peebles, and E.~Price.
\newblock Collision-based testers are optimal for uniformity and closeness.
\newblock \emph{Chic. J. Theor. Comput. Sci}, 25:\penalty0 1--21, 2019.

\bibitem[Friedman and Rafsky(1979)]{fr}
J.~H. Friedman and L.~C. Rafsky.
\newblock Multivariate generalizations of the wald-wolfowitz and smirnov
  two-sample tests.
\newblock \emph{The Annals of Statistics}, pages 697--717, 1979.

\bibitem[Galbraith and Holmes(2012)]{galbraith2012non}
S.~D. Galbraith and M.~Holmes.
\newblock A non-uniform birthday problem with applications to discrete
  logarithms.
\newblock \emph{Discrete Applied Mathematics}, 160\penalty0 (10-11):\penalty0
  1547--1560, 2012.

\bibitem[Henze(1988)]{henzenn}
N.~Henze.
\newblock A multivariate two-sample test based on the number of nearest
  neighbor type coincidences.
\newblock \emph{The Annals of Statistics}, 16\penalty0 (2):\penalty0 772--783,
  1988.

\bibitem[It{\^o}(1951)]{stochasticintegral}
K.~It{\^o}.
\newblock Multiple wiener integral.
\newblock \emph{Journal of the Mathematical Society of Japan}, 3\penalty0
  (1):\penalty0 157--169, 1951.

\bibitem[Janson(1997)]{SJIII}
S.~Janson.
\newblock \emph{Gaussian {H}ilbert spaces}.
\newblock Cambridge University Press, 1997.

\bibitem[Kiani et~al.(2021)Kiani, Gomez-Cabrero, and Bianconi]{networksbiology}
N.~A. Kiani, D.~Gomez-Cabrero, and G.~Bianconi.
\newblock \emph{Networks of Networks in Biology: Concepts, Tools and
  Applications}.
\newblock Cambridge University Press, 2021.

\bibitem[Kim et~al.(2010)Kim, Montenegro, Peres, and Tetali]{discretelogarithm}
J.~H. Kim, R.~Montenegro, Y.~Peres, and P.~Tetali.
\newblock A birthday paradox for markov chains with an optimal bound for
  collision in the pollard rho algorithm for discrete logarithm.
\newblock \emph{The Annals of Applied Probability}, 20\penalty0 (2):\penalty0
  495, 2010.

\bibitem[Kivel{\"a} et~al.(2014)Kivel{\"a}, Arenas, Barthelemy, Gleeson,
  Moreno, and Porter]{kivela2014multilayer}
M.~Kivel{\"a}, A.~Arenas, M.~Barthelemy, J.~P. Gleeson, Y.~Moreno, and M.~A.
  Porter.
\newblock Multilayer networks.
\newblock \emph{Journal of complex networks}, 2\penalty0 (3):\penalty0
  203--271, 2014.

\bibitem[Kwan et~al.(2023)Kwan, Sah, Sauermann, and
  Sawhney]{anitconcentrationramsey}
M.~Kwan, A.~Sah, L.~Sauermann, and M.~Sawhney.
\newblock Anticoncentration in {R}amsey graphs and a proof of the
  {E}rd{\H{o}}s--{M}c{K}ay conjecture.
\newblock \emph{Forum of Mathematics, Pi}, 11:\penalty0 e21, 2023.

\bibitem[Lei et~al.(2020)Lei, Chen, and Lynch]{lei2020consistent}
J.~Lei, K.~Chen, and B.~Lynch.
\newblock Consistent community detection in multi-layer network data.
\newblock \emph{Biometrika}, 107\penalty0 (1):\penalty0 61--73, 2020.

\bibitem[Lov{\'a}sz(2012)]{lovasz_book}
L.~Lov{\'a}sz.
\newblock \emph{Large networks and graph limits}, volume~60.
\newblock American Mathematical Soc., 2012.

\bibitem[Lunag{\'o}mez et~al.(2021)Lunag{\'o}mez, Olhede, and
  Wolfe]{lunagomez2021modeling}
S.~Lunag{\'o}mez, S.~C. Olhede, and P.~J. Wolfe.
\newblock Modeling network populations via graph distances.
\newblock \emph{Journal of the American Statistical Association}, 116\penalty0
  (536):\penalty0 2023--2040, 2021.

\bibitem[MacDonald et~al.(2022)MacDonald, Levina, and Zhu]{macdonald2022latent}
P.~W. MacDonald, E.~Levina, and J.~Zhu.
\newblock Latent space models for multiplex networks with shared structure.
\newblock \emph{Biometrika}, 109\penalty0 (3):\penalty0 683--706, 2022.

\bibitem[Mani and Mikulincer(2024)]{subgraph2024characterizing}
N.~Mani and D.~Mikulincer.
\newblock Characterizing the fourth-moment phenomenon of monochromatic subgraph
  counts via influences.
\newblock \emph{arXiv:2403.14068}, 2024.

\bibitem[Mao et~al.(2021)Mao, Wu, Xu, and Yu]{mao2021testing}
C.~Mao, Y.~Wu, J.~Xu, and S.~H. Yu.
\newblock Testing network correlation efficiently via counting trees.
\newblock \emph{arXiv preprint arXiv:2110.11816}, 2021.

\bibitem[Mao et~al.(2023)Mao, Wu, Xu, and Yu]{mao2023random}
C.~Mao, Y.~Wu, J.~Xu, and S.~H. Yu.
\newblock Random graph matching at otter's threshold via counting chandeliers.
\newblock In \emph{Proceedings of the 55th Annual ACM Symposium on Theory of
  Computing}, pages 1345--1356, 2023.

\bibitem[Mckinney(1966)]{generalizedbirthdayproblem}
E.~H. Mckinney.
\newblock Generalized birthday problem.
\newblock \emph{The American Mathematical Monthly}, 73\penalty0 (4):\penalty0
  385--387, 1966.

\bibitem[Mitzenmacher and Upfal(2017)]{mitzenmacher2017probability}
M.~Mitzenmacher and E.~Upfal.
\newblock \emph{Probability and computing: Randomization and probabilistic
  techniques in algorithms and data analysis}.
\newblock Cambridge university press, 2017.

\bibitem[Mossel et~al.(2010)Mossel, O'Donnell, and Oleszkiewicz]{stability}
E.~Mossel, R.~O'Donnell, and K.~Oleszkiewicz.
\newblock Noise stability of functions with low influences: Invariance and
  optimality.
\newblock \emph{Annals of Mathematics}, pages 295--341, 2010.

\bibitem[Nandi and Stinson(2007)]{nandi2007multicollision}
M.~Nandi and D.~R. Stinson.
\newblock Multicollision attacks on some generalized sequential hash functions.
\newblock \emph{IEEE transactions on Information Theory}, 53\penalty0
  (2):\penalty0 759--767, 2007.

\bibitem[Nourdin and Peccati(2009)]{nourdin2009stein}
I.~Nourdin and G.~Peccati.
\newblock Stein's method on {W}iener chaos.
\newblock \emph{Probability Theory and Related Fields}, 145:\penalty0 75--118,
  2009.

\bibitem[Nourdin and Peccati(2012)]{nourdin2012normal}
I.~Nourdin and G.~Peccati.
\newblock \emph{Normal approximations with {M}alliavin calculus: from {S}tein's
  method to universality}, volume 192.
\newblock Cambridge University Press, 2012.

\bibitem[Nourdin and Rosi{\'n}ski(2014)]{nourdin2014asymptotic}
I.~Nourdin and J.~Rosi{\'n}ski.
\newblock Asymptotic independence of multiple {W}iener-it{\^o} integrals and
  the resulting limit laws.
\newblock \emph{The Annals of Probability}, 42\penalty0 (2):\penalty0 497--526,
  2014.

\bibitem[Nourdin et~al.(2008)Nourdin, Peccati, and Reinert]{rademacher}
I.~Nourdin, G.~Peccati, and G.~Reinert.
\newblock Stein's method and stochastic analysis of {R}ademacher functionals.
\newblock \emph{Electronic Journal of Probability}, 15, 2008.

\bibitem[Nourdin et~al.(2010)Nourdin, Peccati, and
  Reinert]{nourdin2010invariance}
I.~Nourdin, G.~Peccati, and G.~Reinert.
\newblock Invariance principles for homogeneous sums: universality of
  {G}aussian {W}iener chaos.
\newblock \emph{The Annals of Probability}, 38\penalty0 (5):\penalty0
  1947--1985, 2010.

\bibitem[Nualart and Peccati(2005)]{nualart2005central}
D.~Nualart and G.~Peccati.
\newblock Central limit theorems for sequences of multiple stochastic
  integrals.
\newblock \emph{Annals of Probability}, pages 177--193, 2005.

\bibitem[Paninski(2008)]{paninski2008coincidence}
L.~Paninski.
\newblock A coincidence-based test for uniformity given very sparsely sampled
  discrete data.
\newblock \emph{IEEE Transactions on Information Theory}, 54\penalty0
  (10):\penalty0 4750--4755, 2008.

\bibitem[Pedarsani and Grossglauser(2011)]{pedarsani2011privacy}
P.~Pedarsani and M.~Grossglauser.
\newblock On the privacy of anonymized networks.
\newblock In \emph{Proceedings of the 17th ACM SIGKDD international conference
  on Knowledge discovery and data mining}, pages 1235--1243, 2011.

\bibitem[Rosenbaum(2005)]{rosenbaum}
P.~R. Rosenbaum.
\newblock An exact distribution-free test comparing two multivariate
  distributions based on adjacency.
\newblock \emph{Journal of the Royal Statistical Society Series B: Statistical
  Methodology}, 67\penalty0 (4):\penalty0 515--530, 2005.

\bibitem[Rosi\'{n}ski and Samorodnitsky(1999)]{rosinski1999product}
J.~Rosi\'{n}ski and G.~Samorodnitsky.
\newblock Product formula, tails and independence of multiple stable integrals.
\newblock \emph{Advances in stochastic inequalities (Atlanta, GA, 1997)},
  234:\penalty0 169--194, 1999.

\bibitem[Rotar(1979)]{rotar2}
V.~I. Rotar.
\newblock Limit theorems for polylinear forms.
\newblock \emph{Journal of Multivariate analysis}, 9\penalty0 (4):\penalty0
  511--530, 1979.

\bibitem[Schilling(1986)]{schilling}
M.~F. Schilling.
\newblock Multivariate two-sample tests based on nearest neighbors.
\newblock \emph{Journal of the American Statistical Association}, 81\penalty0
  (395):\penalty0 799--806, 1986.

\bibitem[Skeja and Olhede(2024)]{skeja2024quantifying}
A.~Skeja and S.~C. Olhede.
\newblock Quantifying multivariate graph dependencies: Theory and estimation
  for multiplex graphs.
\newblock \emph{arXiv preprint arXiv:2405.14482}, 2024.

\bibitem[Stinson(2005)]{stinson2005cryptography}
D.~R. Stinson.
\newblock \emph{Cryptography: theory and practice}.
\newblock Chapman and Hall/CRC, 2005.

\bibitem[Suzuki et~al.(2006)Suzuki, Tonien, Kurosawa, and
  Toyota]{birthdaymulticollison}
K.~Suzuki, D.~Tonien, K.~Kurosawa, and K.~Toyota.
\newblock Birthday paradox for multi-collisions.
\newblock In M.~S. Rhee and B.~Lee, editors, \emph{Information Security and
  Cryptology -- ICISC 2006}, pages 29--40, Berlin, Heidelberg, 2006. Springer
  Berlin Heidelberg.

\bibitem[\"{U}st\"{u}nel and Zakai(1989)]{ustunel1989independence}
A.~S. \"{U}st\"{u}nel and M.~Zakai.
\newblock On independence and conditioning on {W}iener space.
\newblock \emph{The Annals of Probability}, pages 1441--1453, 1989.

\end{thebibliography}

\normalsize

\appendix

\section{Continuity of the 2-Point Conditional Kernel}

Recall the definition of the 2-point conditional kernel from \eqref{eq:WH}. Note that $W_H$ is symmetric and bounded, but can take values greater 1, hence, $W_H \in \sW_1$, where $\sW_1$ is the space of  bounded, symmetric, measurable functions from $[0, 1]^2$ to $[0, \infty)$. The topology generated by the cut-distance \eqref{eq:Wconvergence} (extended naturally to the space $\sW_1$) is the same as the topology generated by the following norm (see \cite[Lemma 8.11]{lovasz_book}):
\begin{align*}
	||W_1-W_2||_{1 \rightarrow \infty}:=\sup_{f, g: [0, 1] \rightarrow [-1, 1]}\left|\int_{[0, 1]^2} \left(W_1(x, y)-W_2(x, y)\right) f(x) g(y) \mathrm dx \mathrm dy \right| ,
\end{align*}
for $W_1, W_2 \in \sW_1$. The next lemma shows that the map $W \rightarrow W_H$ is Lipschitz in the $\| \cdot \|_{1 \rightarrow \infty}$ norm and, hence, the cut-distance $\| \cdot \|_{\square}$. This, in particular, means that if $G_n$ converges to a graphon $W$, then $W^{G_n}_{H}$ converges in cut-norm to $W_H$.

\begin{lem}\label{lem:W_H_lipschitz}
	The mapping $\eta_H:W\mapsto W_H$ defined on $\sW_1$ is Lipschitz under $\|\cdot\|_{\infty\ra 1}$ and, hence, under $\|\cdot\|_{\square}$.
\end{lem}
\begin{proof} Let $\{U_r\}_{r\in V(H)}$ be a collection i.i.d. $\mathrm{Unif}[0,1]$ random variables. To show the map $\eta_H$ is Lipschitz, it suffices to prove
	\begin{align}\label{eq:WRHsquare}
		\| (W+R)_H- R_H \|_{\infty \rightarrow 1} \lesssim_H \| R \|_{\infty \rightarrow 1} ,
	\end{align}
	for $W, R \in \sW_1$. To this end, fix $1 \leq a \ne b \leq |V(H)|$ and note that
	\begin{align*}
		\Delta_{a, b}(x, y) & := \E\left( \prod_{(i, j)\in E(H)}\left(W(U_i, U_j)+ R(U_i, U_j)\right)-\prod_{(i, j)\in E(H)}W(U_i, U_j) \Bigg| U_a=x, U_b=y \right) \\
		                    & =\E\left(\left. \sum_{\substack{\cA\subseteq  E(H)                                                                                    \\|\cA|\geq 1}} \prod_{(i, j) \in  \cA^c}W(U_i, U_j)\prod_{(i, j) \in  \cA}R(U_i, U_j)\,\right\rvert\, U_a=x, U_b=y \right) .
	\end{align*}
	Recalling \eqref{eq:WH}, note that,
	\begin{align}\label{eq:Deltaxy}
		\| (W+R)_H- R_H \|_{\infty \rightarrow 1}  = \left\| \sum_{1 \leq a \ne b \leq |V(H)|}  \Delta_{a, b} \right\|_{\infty \rightarrow 1} & \leq \sum_{1 \leq a \ne b \leq |V(H)|} \left\| \Delta_{a, b} \right\|_{\infty \rightarrow 1} ,
	\end{align}
	by the triangle inequality. Hence, to prove \eqref{eq:WRHsquare} it suffices to consider each term in the above sum separately. Towards this,
	\begin{align*}
		 & \left\| \Delta_{a, b} \right\|_{\infty \rightarrow 1} \\
		 & = \sup_{\substack{\|f\|_{\infty}\leq 1                \\ \|g\|_{\infty}\leq 1}}\left\lvert \E\left(\E\left(\left. \sum_{\substack{\cA\subseteq  E(H)\\|\cA|\geq 1}} \prod_{(i, j) \in  \cA^c}W(U_i, U_j)\prod_{(i, j) \in  \cA}R(U_i, U_j)\,\right\rvert\, U_a, U_b\right)f(U_a)g(U_b)\right)\right\rvert\\
		 & =\sup_{\substack{\|f\|_{\infty}\leq 1                 \\ \|g\|_{\infty}\leq 1}}\left\lvert \E\left(\sum_{\substack{\cA\subseteq  E(H)\\|\cA|\geq 1}} f(U_a)g(U_b)\prod_{(i, j) \in  \cA^c}W(U_i, U_j)\prod_{(i, j) \in  \cA}R(U_i, U_j)\right)\right\rvert\\
		 & \leq \sum_{\substack{\cA\subseteq  E(H)               \\|\cA|\geq 1}} \sup_{\substack{\|f\|_{\infty}\leq 1\\ \|g\|_{\infty}\leq 1}}\left\lvert   \E\left( f(U_a)g(U_b)\prod_{(i, j) \in  \cA^c}W(U_i, U_j)\prod_{(i, j) \in  \cA}R(U_i, U_j)\right)\right\rvert ,
		\stepcounter{equation}\tag{\theequation}\label{eq:Deltaxyfg}
	\end{align*}
	by the triangle inequality. Now, by a telescoping argument similar to the proof of the counting lemma (see \cite[Theorem 3.7]{graph_limits_I} and \cite[Lemma 10.24]{lovasz_book}) it can be shown that
	\begin{align}\label{eq:WAfg}
		\sup_{\substack{\|f\|_{\infty}\leq 1 \\ \|g\|_{\infty}\leq 1}}\left\lvert   \E\left( f(U_a)g(U_b)\prod_{(i, j) \in  \cA^c}W(U_i, U_j)\prod_{(i, j) \in  \cA}R(U_i, U_j)\right)\right\rvert \lesssim_H  \| R \|_{\infty \rightarrow 1}.\end{align}
	Combining \eqref{eq:Deltaxy}, \eqref{eq:Deltaxyfg}, and \eqref{eq:WAfg}, the result in \eqref{eq:WRHsquare} follows.
\end{proof}

\section{Independence of Bivariate Stochastic Integrals}

Fix $d \geq 1$. Given a bounded function $f: [0, 1]^d \rightarrow \R$, denote by $I_d(f)$ the $d$-dimensional Wiener-It\^{o} stochastic integral of $f$ with respect to a Brownian motion on $[0, 1]$. Conditions under which 2 multiple stochastic integrals are independent are well-known. \citet{ustunel1989independence} provided a useful necessary and sufficient for the independence of multiple stochastic integrals and \citet{rosinski1999product} showed that multiple stochastic integrals
are independent if and only if their squares are uncorrelated.  An asymptotic version of these results was established by \citet{nourdin2014asymptotic}. Using this asymptotic result we prove the following lemma:

\begin{lem}\label{lem:almost_independent_L2}
	Fix $W \in \sW_1$ and consider a sequence of bounded kernels $\{R_n\}_{n \geq 1}$, with $R_n \in \sW_1$, such that $\|R_n\|_{\square} \rightarrow 0$. Then $\bm{Q}_n := (I_2(R_n), I_2(W))$ are asymptotically independent along any subsequence for which $\bm{Q}_n$ has a limit in distribution. \end{lem}
\begin{proof}
	By \cite[Theorem 3.1]{nourdin2014asymptotic}, to show the asymptotic independence it suffices to check the following 2 conditions:
	\begin{align}\label{eq:WR}
		\int_{[0,1]^2} R_n(x, y) W(x, y) \mathrm{d} x \mathrm{d} y \rightarrow 0 \quad \text{ and } \quad \int_{[0, 1]^2} \left(\int_0^1 R_n(x, z) W(z, y) \mathrm{d}z \right)^2 \mathrm{d} x \mathrm{d} y \rightarrow 0.
	\end{align}
	The first condition in \eqref{eq:WR} follows from \cite[Lemma 8.22]{lovasz_book}. For the second condition note that:
	\begin{align}
		 & \int_{[0, 1]^2} \left(\int_0^1 R_n(x, z) W(z, y) \mathrm{d}z \right)^2 \mathrm{d} x \mathrm{d} y \nonumber                                                                                                \\
		 & = \int_{[0, 1]^4} R_n(x, z) W(z, y)  R_n(x, z') W(z', y) \mathrm{d}z \mathrm{d} z' \mathrm{d} x \mathrm{d} y  \nonumber                                                                                   \\
		 & = \int_{[0, 1]^2} \left(\int_{[0, 1]^2} g_n(x, y) R_n(x, z) W(z, y)  \mathrm{d} x\mathrm{d}z\right) \mathrm{d} y  \tag*{(where $g_n(x, y) := \int_{[0, 1]} R_n(x, z') W(z', y) \mathrm{d} z'$)} \nonumber \\
		 & \leq \|R_n\|_{\square} \rightarrow 0. \nonumber
	\end{align}
	This establishes the second condition in \eqref{eq:WR} and completes the proof of Lemma \ref{lem:almost_independent_L2}.
\end{proof}

\section{A CLT for Weighted Sum of Centered $\chi^2$ Random Variables}

In this section we establish a Central Limit Theorem for an infinite weighted sum of independent centered $\chi^2$ random variables.

\begin{lem}\label{lm:chisquareCLT}
	Consider a sequence of infinite sequences $\{(a_s^{(n)}: s\geq 1)\}_{n \geq 1}$ satisfying the following conditions:
	\begin{align}\label{eq:sequence}
		\lim_{n \rightarrow \infty} \max_{s \geq 1} | a_s^{(n)} | = 0   \text{ and } \lim_{n \rightarrow \infty} \sum_{s \geq 1} (a_s^{(n)})^2 = \tau^2 ,
	\end{align}
	for some constant $\tau \geq 0$. Let $\{\xi_s\}_{s \ge 1}$ be an i.i.d. sequence of $\chi^2_k-k$ random variables. Then the infinite sum $\sum_{s \geq 1} a_s^{(n)} \xi_s$ is well defined, for $n$ large enough,  and, as $n \rightarrow \infty$,
	\begin{align}\label{eq:chisquaredsum}
		\sum_{s \geq 1} a_s^{(n)} \xi_s \stackrel{D} \rightarrow N\left(0,2k \tau^2 \right) .
	\end{align}
\end{lem}

\begin{proof}  Fix $M \geq 1$ and let $\cF_M := \sigma(\{\xi_s\}_{s=1}^M)$ be the sigma algebra generated by $\{\xi_s\}_{s=1}^M$ and $S_M^{(n)} := \sum_{s=1}^M a_s^{(n)} \xi_s$. Then $(S_M, \cF_M)$ is a martingale with
	$$\sup_{M \geq 1}\E (S_M^{(n)})^2 \leq  \sum_{s \geq 1} (a_s^{(n)})^2 < \infty,$$
	for all $n$ large enough. Hence, the sum $\sum_{s \geq 1} a_s^{(n)} \xi_s$ is well defined, for $n$ large enough.

	We will show \eqref{eq:chisquaredsum} by establishing the convergence of the Moment Generating Function (MGF). By \cite[Proposition 7.1]{BDM} we know that the MGF of $\sum_{s \geq 1} a_s^{(n)}\xi_i$ is well defined in a neighborhood of zero. In particular, for $|\lambda| < \frac{1}{8}$, one has
	\begin{align*}
		\mathcal{E}_n(\lambda) := \E\left( e^{\lambda \sum_{s \geq 1} a_s^{(n)}\xi_i }  \right) = \prod_{s \geq 1} \frac{e^{- \lambda k a_s^{(n)}}}{\left(1-2 \lambda a_s^{(n)}\right)^{\frac{k}{2}}} .
	\end{align*}
	Taking logarithms on both sides,
	\begin{align*}
		\log \mathcal{E}_n(\lambda) = -\frac{k}{2} \sum_{s \geq 1} \log\left(1-2 \lambda a_s^{(n)}\right)- \lambda k \sum_{s \geq 1} a_s^{(n)}
		 & = \frac{k}{2}\sum_{s \geq 1} \sum_{t=1}^\infty\frac{(2 \lambda a_s^{(n)})^t}{t} - \lambda k \sum_{s \geq 1} a_s^{(n)}        \\
		 & = \frac{k}{2}\sum_{s \geq 1} \sum_{t=2}^\infty\frac{ (2t a_s^{(n)})^t}{t}                                                    \\
		 & = \frac{k}{2}\sum_{s \geq 1} \sum_{t=3}^\infty\frac{(2 \lambda a_s^{(n)})^t}{t} +k \lambda^2 \sum_{s \geq 1} (a_s^{(n)})^2 .
		\stepcounter{equation}\tag{\theequation}
		\label{eq:mgf_doublesum}
	\end{align*}
	Denote by $a^{(n)}_{(s)}$ the $s$-th largest among $\{a_s^{(n)}\}_{s \geq 1}$. Then, for any $L \geq 1$, we have
	\begin{align*}
		4 \tau^2 \geq \sum_{s\geq 1} (a_s^{(n)})^2 \geq\sum_{s=1}^L(a^{(n)}_{(s)})^2\geq L (a^{(n)}_{(L)})^2 ,
	\end{align*}
	for all  $n$ large enough. (Note that the first inequality in the display above holds for all  $n$ large enough by the second condition in \eqref{eq:sequence}.) This implies, $a^{(n)}_{(L)}\leq \frac{2\tau}{\sqrt{L}}$, for $n$ large enough and $L \geq 1$.
	Hence, for $|\lambda| < \frac{1}{8 \tau}$,
	\begin{align}\label{eq:ansum}
		\sum_{s \geq 1} \sum_{t=3}^\infty\frac{|\lambda a_s^{(n)}|^t}{t}  \leq \sum_{s \geq 1} \sum_{t=3}^\infty | \lambda a^{(n)}_{(s)}|^t \leq \sum_{s \geq 1} \sum_{t=3}^\infty\left\lvert\frac{4 \lambda \tau}{\sqrt{s}}\right\rvert^t & \leq \sum_{t=3}^\infty\sum_{s \geq 1} \left\lvert\frac{1}{2\sqrt{s}}\right\rvert^t < \infty .
	\end{align}
	This shows that the double sum in \eqref{eq:mgf_doublesum} is absolutely convergent. Also, for $t \geq 3$,
	\begin{align*}
		\lim_{n \rightarrow \infty} \left|\sum_{s \geq 1} (a_s^{(n)})^t\right| \leq \lim_{n \rightarrow \infty} \lp\max_{s\geq 1} | a_s^{(n)}|\rp^{t-2} \sum_{s \geq 1}  (a_s^{(n)})^2 = 0 ,
	\end{align*}
	from \eqref{eq:sequence}. Hence, by the Dominated Convergence Theorem on the sequence of functions $f_n(t) = \sum_{s\geq 1}|a^{(n)}_{(s)}|^{t}$ (observe that the bound in \eqref{eq:ansum} is uniform over $n$),
	$$ \lim_{n \rightarrow \infty} \sum_{s \geq 1} \sum_{t=3}^\infty\frac{|\lambda a_s^{(n)}|^t}{t} = 0.$$
	Hence, taking limits as $n \rightarrow \infty$ on both sides of \eqref{eq:mgf_doublesum} gives, for $|\lambda| < \frac{1}{8 \tau}$,
	\begin{align*}
		\lim_{n\ra\infty}\mathcal{E}_n(\lambda)= e^{k \tau^2 \lambda^2 } ,
	\end{align*}
	which is the MGF of $N(0, 2k \tau^2)$. This completes the proof of Lemma \ref{lm:chisquareCLT}.
\end{proof}

\end{document}